\documentclass[11pt,a4paper,reqno]{amsart}
\usepackage{tikz}
\usepackage[utf8]{inputenc}
\usepackage{amsmath,bm,bbm,amsthm, amssymb}
\usepackage{fullpage}
\usepackage{color}
\usepackage{animate}
\usepackage{etoolbox}
\usepackage{ulem}
\usepackage{comment}
\usepackage{pgf}
\usepackage{hyperref}
\usetikzlibrary{calc}
\usetikzlibrary{patterns}
\usetikzlibrary{arrows}
\usetikzlibrary{decorations.pathreplacing}
\usepackage[utf8]{inputenc}
\usepackage{dsfont}
\usepackage{float}
\usepackage{enumitem}

\usepackage{pgfplots}
\usepgfplotslibrary{colormaps}
\pgfplotsset{compat=1.17}
\pgfplotsset{
colormap={whiteblack}{gray(0cm)=(1); gray(1cm)=(0)}
}
\usepackage[foot]{amsaddr}

\def\ms{\mathsf}
\def\mc{\mathcal}

\def\mbb{\mathbb}

\def\Z{\mbb Z}
\def\E{\mbb E}

\def\P{\mbb P}
\def\Q{\mbb Q}
\def\R{\mbb R}
\def\N{\mbb N}

\def\PP{\mathcal P}
\def\NN{\mathcal N}
\def\NNN{\mathbf N_\ms{fin}}
\def\NNNN{\mathbf N}

\def\d{{\rm d}}

\def\a{\alpha}

\def\de{\delta}

\def\e{\varepsilon}
\def\eps{\e}

\def\one{\mathbbmss{1}}

\def\bel{\begin{lemma}}
\def\enl{\end{lemma}}
\def\bepr{\begin{proposition}}
\def\enpr{\end{proposition}}
\def\bep{\begin{proof}}
\def\enp{\end{proof}}

\newcommand{\pimsta}[1]{\pi_{-s, t, a}}

\def\to{\uparrow}

\def\th{\theta}

\def\been{\begin{enumerate}}
\def\enen{\end{enumerate}}
\def\im{\item}

	\def\bec{\begin{corollary}}
	\def\enc{\end{corollary}}

\def\bet{\begin{theorem}}
\def\ent{\end{theorem}}

\DeclareMathOperator*{\diam}{diam}
\DeclareMathOperator*{\dist}{dist}

\def\CCC{\ms{CC}}
\def\th{\ms{cr}}
\def\sp{\ms{sp}}
\def\de{\ms{de}}

\theoremstyle{plain}
\newtheorem{theorem}{Theorem}
\newtheorem{proposition}[theorem]{Proposition}
\newtheorem{corollary}[theorem]{Corollary}
\newtheorem{lemma}[theorem]{Lemma}

\theoremstyle{definition}

\theoremstyle{remark}
\newtheorem{remark}{Remark}


\title{Lower large deviations for geometric functionals in sparse, critical and dense regimes}
\date{\today}

\author{Christian Hirsch$^1$}
\address{$^1$Aarhus university}
\email{hirsch@math.au.dk}

\author{Daniel Willhalm$^2$}
\address{$^2$University of Groningen}
\email{d.willhalm@rug.nl}

\begin{document}

\begin{abstract} 
We prove lower large deviations for geometric functionals in sparse, critical and dense regimes. Our results are tailored for functionals with nonexisting exponential moments, for which standard large deviation theory is not applicable. The primary tool of the proofs is a sprinkling technique that, adapted to the considered functionals, ensures a certain boundedness. This substantially generalizes previous approaches to tackle lower tails with sprinkling. Applications include subgraph counts, persistent Betti numbers and edge lengths based on a sparse random geometric graph, power-weighted edge lengths of a $k$-nearest neighbor graph as well as power-weighted spherical contact distances in a critical regime and volumes of $k$-nearest neighbor balls in a dense regime.
\end{abstract}

\maketitle

\section{Introduction}
The theory of large deviations is a central research topic in probability theory which aims to quantify and understand large fluctuations in systems affected by randomness. As it becomes increasingly important to understand the behavior of random systems not only in typical situations but also in unlikely scenarios, large deviations theory has become a central element  in a broad range of application domains, such  as telecommunications, rare-event simulations, insurance mathematics and information theory \cite{dembozeitouni}. While classical large deviations theory predominantly investigates sequences of random variables or time-varying processes, more recently there has been vigorous research activity in investigating large deviations properties of random geometric and topological structures \cite{SchreiberYukich}.

One of the key characteristics of these spatial systems is that we frequently observe a distinctively different behavior in the lower and in the upper large deviation tails. More precisely, for upper large deviations, we often observe \emph{condensation}. That is, the rare events are caused by a highly pathological structure localized in a small part of the sampling window, while the rest of the system behaves essentially as in the typical regime \cite{chatterjeeharel,hirschwill,kerriou}. In contrast, in the lower large deviations, we are typically in a homogenization phase. That means the large deviations are caused by consistent changes away from the typical regime throughout the sampling window. 

The classical techniques to deal with large deviations are predominantly designed to deal with situations where the lower and the upper tails are of the same nature \cite{SchreiberYukich,georgii}. Hence, it is often unclear how to apply them in the geometric situations outlined above. On a mathematical level, the reason for this difficulty is the lack of suitable exponential moments. To address these problems, recently \cite{hirsch} proposed a sprinkling method. Loosely speaking, this method is based on the idea that it is often possible to eliminate pathological configurations through a small modification of the underlying Poisson process. On a technical level, this sprinkling is implemented through a carefully devised coupling construction. The benefit of this sprinkling step is that after this modification, the pathological configurations are removed and become amenable to an analysis with classical tools. 

However, while the examples described in \cite{hirsch} provide a first idea of the feasibility of the sprinkling approach, the assumptions that are imposed prevent the method from being applied to a broad class of models. For instance, while the method in \cite{hirsch} can deal with power-weighted edge lengths of $k$-nearest neighbor graphs, the power is restricted to be smaller than the dimension. In particular, it does not yield the lower-tail complement of the upper tail analysis in \cite{hirschwill}. More generally, the approach in \cite{hirsch} only deals with the critical regime, where the number of relevant Poisson points is proportional to the size of the sampling window. However, in the context of topological data analysis also, different regimes characterized by either much sparser or much denser configurations of points gained substantial interest \cite{kahlemeckes,owadathomas}.

In the present paper, we address the shortcomings described above. More precisely:
\been
\im In the critical regime, we describe an extension of the sprinkling approach that allows us to deal with large deviations of distance-based functionals to a high power.
\im In the sparse regime, we describe the lower large deviations of a large class of additive functionals, including persistent Betti numbers.
\im In the particularly challenging dense regime, we are able to deal with the lower large deviations of large power-weighted $k$-nearest neighbor distances.
\enen
On a methodological level, the key contribution of our work is a substantial improvement of the sprinkling construction from \cite{hirsch}. While in that work, the coupling was relatively basic in the sense that it typically was enough to add a sparsely distributed process of sprinkled points homogeneously throughout the window. In the present paper, we describe sprinkling strategies that are far more adapted to the actual pathological configurations. In particular, in the dense regime, we show that it is even possible to implement a desired coupling in a sequential manner where the distribution of the sprinkling in the next step is allowed to depend on the configuration of the sprinkling constructed so far.

The rest of the present paper is structured as follows. Section \ref{section_model} begins with an introduction of the model and an explanation of how to interpret the different regimes and distinguish them. Next, in Sections \ref{section_model_sparse}, \ref{section_model_thermodynamic} and \ref{section_model_dense}, we give a much more detailed view into every regime, the sparse, critical and dense one, respectively. Each of these subsections also contains requirements for the specific regimes that allow a functional to fit within our frameworks for the lower large deviations and in each subsection a theorem is stated. Afterwards, we give a small overview of the literature that our results build on and identify in which way ours differ from and extend these. Sections \ref{section_applications_thermodynamic} and \ref{section_applications_sparse} then consist of examples of functionals that fit within the frameworks of the sparse and critical regimes. Due to the complexity of the dense regime, we restrict ourselves to the case of volumes of large $k$-nearest neighbor balls. The rest of the paper is devoted to the proofs of the three main theorems for each regime. Section \ref{section_proof_thermodynamic} deals with the proof within the critical regime, Section \ref{section_proof_sparse} with the proof within the sparse regime and Section \ref{section_proof_dense} with the proof within the dense regime.

\section{Model}\label{section_model}
For $d\in\N$, let $\PP_n \subseteq [0,1]^d$ be a Poisson point process with intensity $n$. The unit cube is equipped with the torus distance given by
$$|x-y| := \min_{z\in\Z^d} \|x-y+z\|$$
for $x,y\in[0,1]^d$, where $\|\cdot\|$ represents the Euclidean norm in $\R^d$. For $x\in[0,1]^d$ and $r>0$, we express the closed ball of radius $r$ with respect to the Euclidean or toroidal metric by $B_r(x)$. Which metric is meant will be clear from the context, and we use $\kappa_d$ to denote the volume of the $d$-dimensional unit ball. First, we demonstrate how geometric functionals on the vertex set $\PP_n$ are commonly set up and how to categorize them into one of the three regimes. In general, most geometric functionals, such as subgraph counts of a random geometric graph or power-weighted edge length of the $k$-nearest neighbor graph, can be encoded by a functional of the form
\begin{equation}\label{equation_general_functional}
H_n(\PP_n) := \frac1{s_n} \sum_{X\in\PP_n} \xi_n(X,\PP_n),
\end{equation}
where
\begin{equation}\label{equation_general_score}
\xi_n\colon \R^d \times \NNNN \rightarrow [0,\infty]
\end{equation}
represents the score function, i.e., the contribution of each single vertex of a set of nodes to the whole functional, where by $\NNNN$, we denote the space of locally finite subsets of $\R^d$. Since $\PP_n$ almost surely contains only a finite amount of points, in most cases, it will be sufficient to only define the score function on finite subsets of $\R^d$, which we denote by $\NNN$. In some cases, we desire to only consider such configurations on the torus for which we write $\NNN^{(1)} = \{\varphi\in\NNN\colon \varphi\subseteq [0,1]^d\}$. Further, informally expressed, the normalizing factor $s_n$ corresponds to the expected number of nodes in $\PP_n$ that admit a positive score. We call such points \emph{relevant}. Throughout the paper, we will use the expression $\varphi(A)$ for a configuration $\varphi\in\NNNN$ and a measurable set $A\subseteq\R^d$ to denote the number of points of $\varphi$ that are located within $A$.

We distinguish between three regimes, the sparse, the critical, sometimes also called thermodynamic, and the dense regime. From a heuristic point of view, this distinction comes from the typical amount of Poisson points in the range that determines the score of a relevant point. Loosely speaking, for many score functions, the score of vertices can be determined locally by only looking at a small neighborhood around the considered point. More precisely, the regimes are distinguished by a sequence $(r_n)_n$ such that for a relevant point $X\in\PP_n$, typically
\begin{equation}\label{equation_heuristics_regimes}
\xi_n(X,\PP_n) = \xi_n(X,\PP_n\cap B_{r_n}(X)).
\end{equation}
The simplest case are functionals that represent features of the random geometric graph, in which $r_n$ corresponds to (the order of) the connectivity radius. For this specific example, the asymptotic behavior of the expected degree of a vertex in the random geometric graph characterizes the respective regime. We emphasize that for other functionals, the distinction into the regimes can be more complicated and refer to Sections \ref{section_model_sparse}, \ref{section_model_thermodynamic} and \ref{section_model_dense} for more details about the particular regimes. Sticking with the heuristic explanation and \eqref{equation_heuristics_regimes}, the expected number of Poisson points within the typical range of the score function is consequently of order $n r_n^d$. Hence, there are three possible scenarios for the asymptotics.
\begin{enumerate}
\item Sparse regime: $n r_n^d\overset{n\to\infty}{\longrightarrow} 0$;
\item Critical regime: $n r_n^d\overset{n\to\infty}{\longrightarrow} c >0$;
\item Dense regime: $n r_n^d\overset{n\to\infty}{\longrightarrow} \infty$.
\end{enumerate}

\begin{figure}[H]
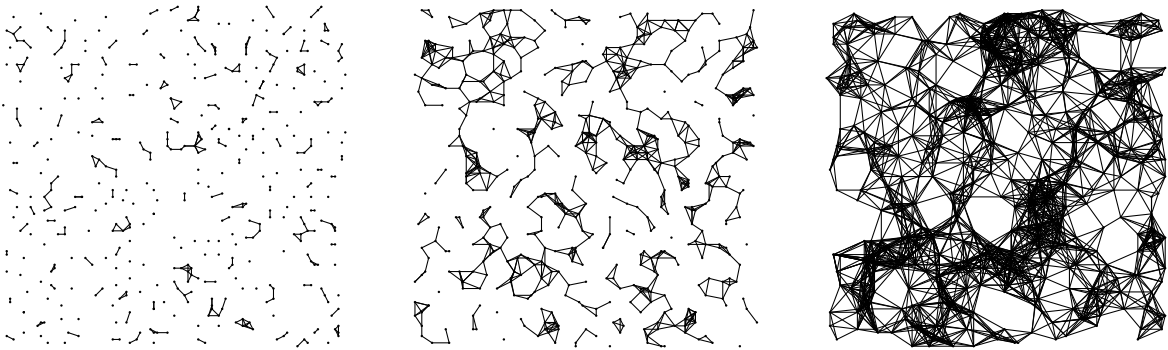

\centering
\input{Tikz_pictures/GeometricGraphSparse}
\hspace{.5cm}
\input{Tikz_pictures/GeometricGraphThermodynamic}
\hspace{.5cm}
\input{Tikz_pictures/GeometricGraphDense} 
\caption{Illustrations of a random geometric graph in a sparse, critical and dense regime.} \label{fig:simulations} 
\end{figure}

The next sections give details about our results in the three different regimes.

\subsection{Sparse regime}\label{section_model_sparse}
In the sparse regime, we investigate functionals for the random geometric graph. We study score functions given by
$$\xi\colon \NNN \rightarrow [0,\infty)$$
defined on finite point configurations in $\R^d$. We also set
\begin{equation}\label{equation_minimum_non_zero}
k_0 := \inf\{m\ge 0\colon \xi(\varphi) >0 \text{ for some } \varphi\subseteq\R^d \text{ with }\#\varphi = m\}
\end{equation}
as the smallest size of a configuration that can yield a positive functional value. We are going to plug configurations of $\PP_n$ into the functional that are rescaled using a sequence of connectivity radii $(r_n)_n\subseteq(0,\infty)$ that will tend to zero. Configurations that have vertices close to the boundary of the torus, which we denote by $\partial[0,1]^d$, might lead to ambiguities if plugged into $\xi$ because the functional itself is not allowed to depend on $n$ and therefore, carries no information about the size of the underlying rescaled torus. For this reason, we generalize the functional to some extent and for $n\in\N$, let
$$\xi_n\colon\NNN^{(1)} \rightarrow [0,\infty)$$
be a functional such that for all configurations $\varphi\in\NNN^{(1)}$ with $\dist(\varphi,\partial[0,1]^d) > r_n$
$$\xi_n(\varphi) = \xi(r_n^{-1} \varphi),$$
where $\dist(\cdot,\cdot)$ denotes the Euclidean distance between two subsets of $\R^d$.

We require $(\xi_n)_n$ and $\xi$ to satisfy the following conditions that are related to the requirements in \cite[Section 3]{hirschowada}.
\begin{enumerate}
\item $\xi$ is translation invariant. That means for all $\varphi\in\NNN$ and $y\in\R^d$, a shift of the configuration $\varphi$ with the vector $y$ does not affect its value, i.e,
\begin{equation}\label{INV}\tag{\textbf{INV}}
\xi(\varphi+y)=\xi(\varphi).
\end{equation}
\item $\xi$ is locally determined for configurations of size $k_0$, which means for all $\varphi\in\NNN$ with $\#\varphi=k_0$
\begin{equation}\label{LOC}\tag{\textbf{LOC}}
\xi(\varphi) = 0 \qquad\text{if } \diam(\varphi) > k_0,
\end{equation}
where $\diam(\varphi) := \max_{y\neq z\in \varphi} \|y-z\|$ denotes the maximal Euclidean distance between points in $\varphi$.
\item For each $m>0$ there exists $b:=b(m)>0$ such that for every $n\in\N$ and every configuration $\varphi\in\NNN^{(1)}$, it holds that
\begin{equation}\label{BND}\tag{\textbf{BND}}
\xi_n(\varphi)\le b \qquad\text{when } \#\varphi \le m.
\end{equation}
\item It holds that
\begin{equation}\label{POS}\tag{\textbf{POS}}
\int_{\R^{(k_0-1)d}} \xi(\{0,x_2,\dots,x_{k_0}\}) \d(x_2,\dots,x_{k_0}) > 0.
\end{equation}
\end{enumerate}

Requiring \eqref{INV} does not exclude any common functionals that represent statistics of random geometric graphs. \eqref{LOC} can be interpreted as a condition that validates $k_0$ as smallest size of a connected component with positive score and is implied if $\xi$ is additive, see Remark \ref{remark_sparse_functional}. Condition \eqref{BND} yields that the score of finite-sized components is finite and \eqref{POS} is a technical condition needed for the result in \cite[Section 3]{hirschowada} that we are going to invoke.

Henceforth, $\ms{GG}_n(\varphi)$ denotes the geometric graph with respect to $|\cdot|$ and with connectivity radius $r_n$ on $\varphi\in\NNN^{(1)}$. Now, we consider the lower large deviations of the functional
\begin{equation}\label{equation_functional_RGG}
H_n^\sp := H_n^\sp(\PP_n) := \frac{1}{\rho_{n,k_0}^\sp} \sum_{\varphi\subseteq \PP_n} \xi_n(\varphi) s_n(\varphi,\PP_n).
\end{equation}
Here, $s_n(\varphi,\PP_n)$ is an indicator function, taking value $1$ if $\varphi$ is a connected component $\ms{GG}_n(\PP_n)$, i.e., for $\varphi\subseteq\psi\in\NNN^{(1)}$ that indicator is given by
\begin{equation}\label{equation_indicator_torus}
s_n(\varphi,\psi) := \one\{\varphi\text{ is a connected component of }\ms{GG}_n(\psi)\}.
\end{equation}
The configuration $r_n^{-1} \varphi$ for $\varphi\subseteq\PP_n$ is considered as a subset of the torus $[0,r_n^{-1}]^d / \sim$ and the normalizing factor has the form
$$\rho_{n,k_0}^\sp := n^{k_0} r_n^{d(k_0-1)},$$
which can be interpreted as the order of the expected number of points that are part of some connected component of size $k_0$.

\begin{remark}\label{remark_sparse_functional}
\begin{enumerate}
\item Note that the functional in \eqref{equation_functional_RGG} is stated in a more general form than suggested in \eqref{equation_general_functional} and \eqref{equation_general_score}. We could recover the representation that sums over all nodes of $\PP_n$ with an indicator that is only nonzero for one vertex of each connected component.
\item Most examples of such functionals, such as subgraph counts, Betti numbers and edge lengths, also fulfill that $\xi_n$ is additive for all $n\in\N$, which means 
$$\xi_n(\varphi_1 \cup \varphi_2) = \xi_n(\varphi_1) + \xi_n(\varphi_2)$$
whenever the distance between $\varphi_1\in\NNN^{(1)}$ and $\varphi_2\in\NNN^{(1)}$ with respect to the toroidal metric is larger than $r_n$. For such functionals we could also write the functional $H_n^\sp$ as $\xi_n(\PP_n)/\rho_{n,k_0}^\sp$.
\end{enumerate}
\end{remark}

Next, along the lines of \cite{hirschowada}, we define a measure on the set $(0,\infty)$ by
$$\tau_{k_0}^\sp(A) := \frac1{k_0!} \lambda_{k_0-1}(\{(y_2,\dots,y_{k_0})\in \R^{d(k_0-1)} \colon \xi(\{0,y_2,\dots,y_{k_0}\})\in A\}),$$
for a measurable $A\subseteq (0,\infty)$, where $\lambda_{k_0-1}$ corresponds to the Lebesgue measure on $\R^{d(k_0-1)}$. Additionally, define the relative entropy of a Radon measure $\rho$ on $(0,\infty)$ by
$$h^\sp(\rho \mid \tau_{k_0}^\sp) := \begin{cases} \int_{(0,\infty)} \log \frac{\d \rho}{\d \tau_{k_0}^\sp}(x) \rho(\d x) - \rho((0,\infty)) + \tau_{k_0}^\sp((0,\infty)) &\text{if }\rho \ll \tau_{k_0}^\sp \\
\infty &\text{otherwise}
\end{cases},$$
where $\rho \ll \tau_{k_0}^\sp$ denotes absolute continuity of $\rho$ with respect to $\tau_{k_0}^\sp$. Note that in accordance with \cite[Remark 3.6]{hirschowada}, under some circumstances, some simplifications of the rate function are possible. We refer to the examples in Section \ref{section_applications_sparse} for details.

The first main theorem states that $H_n^\sp$ admits lower large deviations with rate function $h^\sp(\cdot \mid \tau_{k_0}^\sp)$.

\begin{theorem}[Lower large deviations in the sparse regime]\label{theorem_main_sparse}
Assume that \eqref{INV}, \eqref{LOC}, \eqref{BND} and \eqref{POS} are satisfied and assume that $k_0\in[1,\infty)$. If $n r_n^d\rightarrow 0$ and $\rho_{n,k_0}^\sp\rightarrow\infty$, then, for $a\in\R$
\begin{equation}\label{upper_bound_sparse}
\limsup_{n\to\infty} \frac{1}{\rho_{n,k_0}^\sp} \log \P(H_n^\sp \leq a) \leq -\inf_{\rho\colon T^\sp(\rho) \le a} h^\sp(\rho \mid \tau_{k_0}^\sp)
\end{equation}
and
\begin{equation}\label{lower_bound_sparse}
\liminf_{n\to\infty} \frac{1}{\rho_{n,k_0}^\sp} \log \P(H_n^\sp < a) \geq -\inf_{\rho\colon T^\sp(\rho) < a} h^\sp(\rho \mid \tau_{k_0}^\sp).
\end{equation}
where $T^\sp(\rho) := \int_{(0,\infty)} x \d\rho(x)$.
\end{theorem}

If we assume that $n r_n^d \overset{n\to\infty}{\longrightarrow} 0$, we are indeed in a sparse random geometric graph. But, also using our characterization of the regimes, this setting deserves to be labeled sparse. To verify this, we give a small outlook on the proof of the lower large deviations in this case. First, the typical range to determine the score of a node corresponds to the typical size of a connected component. As it turns out, connected components of size $k_0+1$ or larger do not significantly contribute to the lower large deviations. Therefore, typically the range we have to consider to determine the score of a node or rather the volume occupied by a typical component size is bounded by $k_0^d n r_n^d$, which tends to $0$.

\subsection{Critical regime}\label{section_model_thermodynamic}
For the critical regime, we let $\xi$ be a measurable function
$$\xi \colon \R^d\times\NNNN \rightarrow [0,\infty].$$
Its desired properties are specified later. To turn $\xi$ into the score function we scale everything with the factor $n^{1/d}$ and define
\begin{equation}\label{equation_score_thermo}
\xi_n \colon \R^d\times\NNNN \rightarrow [0,\infty],\ (x,\varphi) \mapsto\xi(n^{1/d} x,n^{1/d} \varphi).
\end{equation}
Here, unlike the sparse regime, we give two different forms of the functional of interest.

\noindent\textbf{Representation A:} We can sum up the scores of each node of the Poisson point process, which is encoded by
\begin{subequations}
\begin{equation}\label{equation_functional_thermo1}
H_n^\th := H_n^\th(\PP_n) := \frac1{n} \sum_{X\in\PP_n} \xi_n(X,\PP_n).
\end{equation}
\textbf{Representation B:} It is also possible to integrate the scores of all space points in $[0,1]^d$, which can be represented by
\begin{equation}\label{equation_functional_thermo2}
H_n^\th := H_n^\th(\PP_n) := \int_{[0,1]^d} \xi_n(x,\PP_n) \d x.
\end{equation}
\end{subequations}

Power-weighted edge lengths of $k$-nearest neighbor graphs is an example of a functional that can be displayed using representation \textbf{A}. Spherical contact distances of space points can be encoded with representation \textbf{B}. See, Section \ref{section_applications_thermodynamic} for details.

\begin{remark}
It is possible to express every functional in representation \textbf{A} in terms of representation \textbf{B} and treat \eqref{equation_functional_thermo1} as a special case of \eqref{equation_functional_thermo2} by using that
\begin{equation}\label{equation_sum_integral_rep}
\frac1{n} \sum_{X\in\PP_n} \xi_n(X,\PP_n) = \int_{[0,1]^d} \sum_{y\in\PP_n\cap B_{(n\kappa_d)^{-1/d}}(x)} \xi_n(y,\PP_n) \d x,
\end{equation}
which can be verified by an application of Fubini's theorem. If all our requirements for a score function would directly translate to the sum of the score function over nodes in a small volume, we could solely consider representation \textbf{B}.
However, we aim to study the lower large deviations of functionals for which some of the requirements for the score function do not translate. In particular, the sum over scores of nodes in a small space can be excessively large if there are many nodes, even if the individual scores are bounded. For this reason, we chose to use two different representations.
\end{remark}

In the critical regime, the notion of \emph{stabilization} plays an important role in many frameworks that deal with limit theory for geometric functionals, see, for example, \cite{weakLLN} or \cite{SchreiberYukich}. Namely, let a function
$$\mc R\colon \R^d\times\NNN \mapsto [0,\infty]$$
be homogeneous of degree $1$, which means that for all $m>0$, $\varphi\in\NNN$ and $x\in\varphi$ it holds that
\begin{equation}\label{equation_stabilization_radius_homogeneity}
\mc R(mx, m\varphi) = m \mc R(x, \varphi).
\end{equation}
Further, we ask for events of the form $\{\mc R(x,\PP_n) \le r\}$ to be measurable with respect to $\PP_n\cap B_r(x)$ for each $x\in[0,1]^d$ and $r>0$.
We call $\mc R$ stabilization radius for $\xi$ if for every $n\in\N$ and $x\in[0,1]^d$
\begin{equation}\label{equation_stabilization_radius}
\P\big(\xi_n(x,\PP_n) = \xi_n(x,\PP_n\cap B_{\mc R(x,\PP_n)}(x))\big) = 1.
\end{equation}

To be able to apply sprinkling to couple two Poisson processes, for an $M>0$ and each $n\in\N$, we introduce $\PP_n^{-,M}$ as a thinning of $\PP_n$ with survival probability $1-M^{-1}$, as well as $\PP_n^{+,M}$ as a Poisson point process on $[0,1]^d$ with intensity $n M^{-1}$ that is independent of $\PP_n$ and the thinning. Then,
\begin{equation}\label{equation_union_poisson_prcesses_thermo}
\PP_n^M := \PP_n^{-,M}\cup \PP_n^{+,M},
\end{equation}
is a Poisson point process $\PP_n^M$ on $[0,1]^d$ with the same distribution as $\PP_n$. The goal for the applications will be to let $\PP_n^{-,M}$ fully cover $\PP_n$ and to sprinkle in additional nodes using $\PP_n^{+,M}$ to control the stabilization radii while at the same time $H_n^{M}(\PP_n^M)$ approximates $H_n^{M}(\PP_n)$. For this purpose, we define an event that is supposed to be the goal of the sprinkling. Here, we need to distinguish between the two representations \eqref{equation_functional_thermo1} and \eqref{equation_functional_thermo2} because in the former, only the nodes of the Poisson point process need to stabilize after the sprinkling.

\noindent\textbf{Representation A:} In the first case, we define the event
\begin{subequations}
\begin{equation}\label{equation_stabilization_event1}
E_n^M := \Big\{\sup_{X\in\PP_n} \mc R(X,\PP_n\cup \PP_n^{+,M})\leq M/n^{1/d}\Big\}
\end{equation}
that the maximal stabilization radius of a node of  $\PP_n\cup\PP_n^{+,M}$ is bounded by $M/n^{1/d}$.

\noindent\textbf{Representation B:}
In the second case, we let
\begin{equation}\label{equation_stabilization_event2}
E_n^M := \Big\{\sup_{x\in [0,1]^d} \mc R(x,\PP_n\cup \PP_n^{+,M})\leq M/n^{1/d}\Big\}
\end{equation}
\end{subequations}
be the set that the maximal stabilization radius of a space point in $[0,1]^d$ with respect to $\PP_n\cup\PP_n^{+,M}$ is bounded by $M/n^{1/d}$. We note that here, $E_n^M$ might not be measurable. But this is of no concern because we only have to deal with subsets of $E_n^M$ later that certainly will be measurable.

Next, for a functional $\xi$ to fit in our framework for lower large deviations in the critical regime, we require additional conditions. Condition \eqref{STA} limits the magnitude of a score function conditioned on a bounded stabilization radius. \eqref{INC} makes tools such as monotone convergence available to use in the proof. \eqref{STA} and \eqref{INC} are satisfied by most examples of score functions in the literature. \eqref{SPR1}, \eqref{SPR2} and \eqref{SPR3} are more restrictive. They make sure that it is possible to find a strategy for sprinkling that bounds the maximal stabilization radius without creating too much excess in the functional. Details about the specific strategies are given in Section \ref{section_applications_thermodynamic}.
\begin{enumerate}
\item Let there exist a stabilization radius $\mc R$ for $\xi$ such that for $x\in[0,1]^d$ and $M>0$ large enough and $n\in\N$
\begin{equation}\label{STA}\tag{\textbf{STA}}
\P\big(\mc R(x,\PP_n) \leq M/n^{1/d}, \xi_n(x, \PP_n) > g(M)\big) = 0
\end{equation}
for some function $g\colon(0,\infty)\rightarrow (0,\infty)$. In particular, $\mc R$ has to satisfy \eqref{equation_stabilization_radius_homogeneity} and \eqref{equation_stabilization_radius}.
\item For each $r>0$, there exists a functional $\xi^r\colon\R^d\times\NNNN \rightarrow [0,\infty]$ bounded by some $r$-dependent constant such that for each $\varphi\in\NNNN$ and $x\in\R^d$ it holds that $\xi^r(x,\varphi)=\xi^r(x,\varphi\cap B_r(x))$ and
\begin{equation}\label{INC}\tag{\textbf{INC}}
\xi^r(x,\varphi) \to \xi(x,\varphi),
\end{equation}
as $r\rightarrow\infty$. In words, $\xi^r$ is nondecreasing with pointwise limit $\xi$.
\end{enumerate}

\noindent Before the last set of requirements, for each $n\in\N$ and $M>0$, we introduce two cut-off versions of the score function using the map $g$ from \eqref{STA} by
$$\xi_n^{M', M}(x, \varphi) = \xi(n^{1/d}x, (n^{1/d}\varphi\cap B_{M'}(n^{1/d} x))) \wedge g(M)$$
and $\xi_n^{M}(x,\varphi) := \xi_n^{M, M}(x, \varphi)$
where $x\in\varphi\in\NNN$. Then, for representation \textbf{A}, we write
\begin{subequations}
\begin{equation}
H^{M',M}_n := H_n^{M',M}(\PP_n) := \frac1{n} \sum_{X\in\PP_n} \xi_n^{M',M}(X,\PP_n)
\end{equation}
and for  representation \textbf{B},
\begin{equation}
H^{M',M}_n := H_n^{M',M}(\PP_n) := \int_{[0,1]^d} \xi_n^{M',M}(x,\PP_n) \d x
\end{equation}
\end{subequations}
as well as $H^{M}_n := H^{M,M}_n$ in both cases for the respective functionals.

\begin{enumerate}[resume]
\item Define the event
$$F_n^{M,(1)} := \{\PP_n = \PP_n^{+,M}\}$$
and for a collection of positive integers $m$ and $I_n^M(\PP_n)$, and a family of disjoint balls in $[0,1]^d/\sim$ that may depend on the Poisson point process
$$(B_{n,i}^M(\PP_n))_{i\in\{1,\dots, I_n^M(\PP_n)\}}$$ with volume $V / n$ for some $V>0$, we set
\begin{equation}\label{equation_sprinkling_event1}
F_n^{M,(2)} := \big\{\PP_n^{+,M}\big([0,1]^d \setminus (\cup_{i = 1}^{I_n^M(\PP_n)} B_{n,i}^M(\PP_n))\big) = 0\big\}
\end{equation}
and
\begin{equation}\label{equation_sprinkling_event2}
F_n^{M,(3)} := \bigcap_{i = 1}^{I_n^M(\PP_n)} \{\PP_n^{+,M}(B_{n,i}^M(\PP_n)) = m\}.
\end{equation}
We assume that the functional allows for such a collection such that
\begin{enumerate}
\item for $M$ sufficiently large, we have
\begin{equation}\label{SPR1}\tag{\textbf{SPR1}}
F_n^M := F_n^{M,(1)} \cap F_n^{M,(2)} \cap F_n^{M,(3)}
 \subseteq E_n^M;
\end{equation}
\item under $\{H_n^{M',M}<a\}$, for $a\in\R$, there exists $c_M^{(1)} \in o(1/\log M)$ as $M\rightarrow\infty$ satisfying that almost surely
\begin{equation}\label{SPR2}\tag{\textbf{SPR2}}
I_n^M(\PP_n) \le c_M^{(1)} n;
\end{equation}
\item there exists $c_M^{(2)} \in o(1)$ as $M\rightarrow\infty$ satisfying that under $\{H_n^{M',M}<a\}\cap F_n^M$, for $a\in\R$, almost surely
\begin{equation}\label{SPR3}\tag{\textbf{SPR3}}
H_n^M(\PP_n^M) \le H_n^{M',M}(\PP_n) + c_M^{(2)}
\end{equation}
if $M$ is sufficiently large.
\end{enumerate}
\end{enumerate}

Similar to \cite{hirsch}, we give the rate function in its entropy-based formulation. For a stationary point process $\mc Q$ defined on $\R^d$, we let $\Q$ be its law, $\Q_n$ be the law $\Q$ restricted to the cube $[0,n^{1/d}]^d$ and $\P_n$ be the law of $n^{1/d}\PP_n$. This lets us set
$$h^\th(\Q) := \begin{cases}\lim_{n\to\infty} \frac{1}{n} \int_\NNNN \log \frac{\d\Q_n}{\d\P_n}(\varphi) \d\Q_n(\varphi) &\text{if }\Q_n \ll \P_n \\
\infty &\text{otherwise}
\end{cases}.$$
Further, for any measure $\widetilde\Q$ on $\NNNN$, we use $\widetilde\Q[\xi]$ to denote $\int_\NNNN \xi(0,\varphi) \d\widetilde\Q(\varphi)$. Next, we need to introduce the Palm version of $\Q$. As it is stated in \cite{georgii}, $\Q$ with finite intensity has a unique finite measure on $\NNNN$ that we denote by $\Q^o$, the Palm version, with the property that for all measurable functions $f\colon\R^d\times\NNNN \rightarrow [0,\infty)$ the equation
$$\E_\Q \Big[\sum_{x\in\varphi} f(x,\varphi - x)\Big] = \int_{\R^d} \int_\NNNN f(x,\varphi\cup\{x\}) \d\Q^o(\varphi) \d x$$ is fulfilled.

This lets us state the theorem dealing with the lower large deviations for the critical regime.

\begin{theorem}[Lower large deviations in the critical regime]\label{theorem_main_thermodynamic}
Let $a>0$.
\begin{itemize}
\item[a)] Assume that $\xi$ satisfies \eqref{INC}. Then,
\begin{equation}\label{upper_bound}
\limsup_{n\to\infty} \frac1{n} \log \P(H_n^\th \leq a) \leq -\inf_\Q h^\th(\Q),
\end{equation}
where the infimum expands over $\{\Q\colon\Q^o[\xi]\le a\}$ or $\{\Q\colon\Q[\xi]\le a\}$ for representation \textbf{A} and representation \textbf{B}, respectively.
\item[b)] Let $H_n^\th$ be given either in representation \textbf{A} or representation \textbf{B}. Assume that $\xi$ satisfies \eqref{STA}, \eqref{SPR1}, \eqref{SPR2} and \eqref{SPR3} for the respective form of $H_n^\th$. Then,
\begin{equation}\label{lower_bound}
\liminf_{n\to\infty} \frac1{n} \log \P(H_n^\th < a) \geq -\inf_\Q h^\th(\Q),
\end{equation}
where the infimum expands over $\{\Q\colon\Q^o[\xi]<a\}$ or $\{\Q\colon\Q[\xi]<a\}$ for representation \textbf{A} and representation \textbf{B}, respectively.
\end{itemize}
\end{theorem}

To see that this coincides with our characterization of the critical regime, we first point out that in order to categorize functionals in representation \textbf{B} into a regime, the characterization via relevant nodes needs to be extended. When dealing with an integral instead of a sum it is sensible to consider any space point $x\in[0,1]^d$ in terms of relevance. Assuming the integral representation for now, we recall that the stabilization radius $\mc R$ is homogeneous of order $1$. In particular, for any relevant $x\in[0,1]^d$ we observe that
$$n^{1/d}\mc R(x,\PP_n) = \mc R(n^{1/d}x,n^{1/d}\PP_n).$$
Note that $n^{1/d}\PP_n$ is a Poisson point process on $[0,n^{1/d}]^d$ with intensity $1$, and thus, for large $n$ typically $\mc R(n^{1/d}x,n^{1/d}\PP_n)$ does not depend on $n$ anymore. Thus, typically $\mc R(x,\PP_n)$ should be of order $n^{-1/d}$, and therefore also the typical range that we need to consider to determine a score of a relevant point, which justifies classifying this framework as critical. If we only consider relevant nodes $X\in\PP_n$, we can repeat the same steps for representation \textbf{A}.

Before continuing with the dense case, we briefly elaborate on the representation of the score function in \eqref{equation_score_thermo}. If a score function is homogeneous of degree $\beta$, thus, there exists $\beta\in\R$ such that for all $m>0$ and $\varphi\in\NNN$ and $x\in\varphi$ it holds that $\xi(m x, m \varphi) = m^\beta \xi(x, \varphi)$, then, the rescaling by $n^{1/d}$ in the arguments of the score function could be replaced by a different normalizing factor for the functional. Power-weighted edge lengths of $k$-nearest neighbor graphs are such an example.

\subsection{Dense regime}\label{section_model_dense}
Since the case of dense spatial networks requires much finer technical argumentation, we focus only on one type of functional for a $k$-nearest neighbor graph for an arbitrary $k\in\N$. In particular, we associate the $k$-nearest neighbor graph with the functional representing large volumes of $k$-nearest neighbor balls. For $x\in\varphi\in\NNN$, this is encoded in
\begin{equation}\label{equation_stabilization_dense}
R_k(x, \varphi) := \inf\{r > 0 \colon \varphi( B_r(x) \setminus\{x\}) \ge k\}.
\end{equation}
This lets us define the according functional by
\begin{equation}\label{equation_functional_kNN}
H_n^\de := H_n^\de(\PP_n) :=  \frac{1}{\rho_{n,k}^\de} \sum_{X\in \PP_n} (n\kappa_d R_k(X, \PP_n)^d - a_n - s_0)_+,
\end{equation}
where $s_0\in\R$ and $(a_n)_n\subseteq\R$ is a sequence that tends to infinity slower than $n$. The normalizing factor has the form
$$\rho_{n,k}^\de := n a_n^{k-1} e^{-a_n}.$$
This factor is derived from the computation
$$\P\big(R_k(X, \PP_n)^d \ge a_n/(n\kappa_d)\big) = \sum_{i=0}^{k-1} \frac{a_n^{i-1}}{i!} e^{-a_n}$$
and represents the expected number of points for which the maximum in \eqref{equation_functional_kNN} is nonzero.

We proceed as in \cite{hirschowadakang} and define a measure on $E_0 := [s_0,\infty)$ by
$$\d\tau_k^\de(x) := \frac{e^{-x}}{(k - 1)!} \d x$$
and denote the relative entropy of a Radon measure $\rho$ on $E_0$ with respect to $\tau_k^\de$ by
$$h^\de(\rho \mid \tau_k^\de) = \begin{cases} \int_{E_0} \log \frac{\d \rho}{\d \tau_k^\de}(x) \d\rho(x) - \rho(E_0) + \tau_k^\de(E_0) &\text{if }\rho \ll \tau_k^\de \\
\infty &\text{otherwise}
\end{cases}.$$

This lets us state the lower large deviations for the functional in \eqref{equation_functional_kNN}.
\begin{theorem}[Lower large deviations in the dense regime]\label{theorem_main_dense}
	Let $(a_n)_n$ be a sequence such that $a_n\rightarrow\infty$ and $a_n - \log n - (k - 1) \log\log n \rightarrow -\infty$. Then, for $a\in\R^d$
\begin{equation}\label{upper_bound_dense}
\limsup_{n\to\infty} \frac{1}{\rho_{n,k}^\de} \log \P(H_n^\de \leq a) \leq -\inf_{\rho\colon T_k^\de(\rho) \le a} h^\de(\rho \mid \tau_k^\de)
\end{equation}
and
\begin{equation}\label{lower_bound_dense}
\liminf_{n\to\infty} \frac{1}{\rho_{n,k}^\de} \log \P(H_n^\de < a) \geq -\inf_{\rho\colon T_k^\de(\rho) < a} h^\de(\rho \mid \tau_k^\de),
\end{equation}
	where $T_k^\de(\rho) := \int_{E_0} x-s_0 \d\rho(x)$.
\end{theorem}

We point out that for a node to have a positive score within any configuration, we have to consider a range of at least $r_n :=((a_n+s_0)/(n\kappa_d))^{1/d}$. Then, $nr_n^d$ diverges if $a_n\rightarrow\infty$. Therefore, typically we would expect to consider an infinite amount of points in the volume within range, and thus, calling this regime dense is indeed sensible.

\subsection{Outline}
Lower large deviations or even large deviation principles for geometric functionals have been derived for sparse, critical and dense regimes in \cite{hirschowada}, \cite{hirsch} and \cite{hirschowadakang}. To achieve an extension of those results, we rely on the technique of sprinkling \cite{sprinkling}, which was already successfully used as a main tool to prove lower large deviations in \cite{hirsch}. In general, it means that we carefully perform small changes to the underlying process at locations that we deem as not suitable in a way such that the functional applied to the adapted configuration still approximates the one with the original point configuration. Mathematically speaking, the idea behind it is to couple two Poisson point processes such that conditioned on one of them, applying the functional to the other one guarantees some additional properties of the score function that allow us to invoke general large deviations theory. In the following paragraphs, we give an overview of the extensions of the sprinkling technique derived in the present work compared to the results from \cite{hirsch}, \cite{hirschowada} and \cite{hirschowadakang}.
\begin{enumerate}
\item \emph{Critical regime:} For the critical case, \cite{hirsch} applies sprinkling on a macroscopic level to control the maximal stabilization radius of any node without significantly altering the functional. A coupled Poisson point process retains all nodes from the original process and consistently inserts additional points across the observation window. The results in \cite{hirsch} are limited to certain functionals for which the magnitude of the score function is comparable to the $d$th power of the stabilization radius. For instance, power-weighted edge lengths for a power as large as or larger than $d$ do not meet the requirements for their results. This restriction substantially simplifies the analysis because in that case regularly inserting points does not alter the functional by a big margin. We will examine some functionals that violate this condition, which requires a much finer adaption of the sprinkling to the studied functional as we will demonstrate in Section \ref{section_applications_thermodynamic}.
\item \emph{Sparse regime:} For a sparse random geometric graph, \cite{hirschowada} derives a large deviations principle for empirical measures counting potentially connected components of a fixed size and certain statistics derived from these. Their strategy builds on weak dependencies of scores assigned to relatively distant connected components in the sparse setting. This lets them approximate functionals restricted to each single box with i.i.d.\ Poisson random measures and apply well-established large deviations theory. However, for their proof to work, it is necessary that considered components cannot be too big. Otherwise, the exponential moments cannot be handled anymore. Using sprinkling, we extend their results. The framework that we present in Section \ref{section_model_sparse} for the sparse regime also focuses on functionals for the random geometric graph but allows to consider connected components of arbitrary size.
\item \emph{Dense regime:} For an empirical measure counting large $k$-nearest neighbor distances, \cite{hirschowadakang} provides a large deviation principle. It proceeds similarly to \cite{hirschowada} by introducing a grid and by approximating the restricted functionals. In our extension, presented in Section \ref{section_model_dense}, we aim to leave the empirical measure setting and use \cite[Theorem 2.1]{hirschowadakang} combined with a sprinkling argument to derive lower large deviations for the functional that directly represents the sum of large distances to the $k$-closest point. The general way sprinkling is applied here is similar to the sparse case. However, due to the finer dependencies between adjacent boxes that have to be resampled, the procedure becomes much more complicated. For this reason, we go sequentially through the boxes, deciding whether to resample them and also making sure that each box, if resampled or not, does not affect the potential resampling of the next boxes negatively.
\end{enumerate}
\section{Examples}
\subsection{Functionals for critical spatial random networks}\label{section_applications_thermodynamic}

\subsubsection{Power-weighted edge lengths of the directed k-nearest neighbor graph}\label{section_applications_thermodynamic_kNN}
Let $k\in\N$ and $\a\geq 0$ be arbitrary. In the directed $k$-nearest neighbor graph, there is a directed edge from each node to its $k$ closest neighbors. We aim to represent the statistic of the power-weighted edge lengths using representation \textbf{A}.
To achieve that, for each $n\in\N$, we let the score function $\xi_n$ be given by
$$\xi(x,\varphi) := \sum_{y\in\varphi\cap B_{R_k(x,\varphi)}(x)} \|x-y\|^\a$$
for $x\in\R^d$ and $\varphi\in\NNNN$, where we recall $R_k$ from \eqref{equation_stabilization_dense} in the dense case, which simultaneously acts as stabilization radius here. For formality reasons, we set $\xi(x,\varphi) := \infty$ if $\#\varphi<k$. Note that when we plug $\PP_n$ into the functional $\xi_n$, we replace the Euclidean norm $\|\cdot\|$ with the toroidal distance on $[0,n^{1/d}]^d/\sim$. Further, the case $\alpha<d$ was already covered in \cite{hirsch}. This functional satisfies \eqref{STA} with the choice $g(m):=e^m$ and also \eqref{INC} is satisfied when choosing $\xi^r(x,\varphi) := \xi(x,\varphi\cap B_r(x))\wedge r$ for $r>0$, $x\in\R^d$ and $\varphi\in\NNNN$. 

In order to show that the sprinkling requirements \eqref{SPR1}, \eqref{SPR2} and \eqref{SPR3} hold as well, we denote all nodes with exceptionally large stabilization radii by
$$\mc J_n^M := \mc J_n^M := \{X\in \PP_n \colon \mathcal{R}(X,\PP_n) > M_n\},$$
where we use the abbreviation $M_n := M/n^{1/d}$.
We point out that the number of vertices in $\PP_n$ on the torus $[0,1]^d/\sim$ with a stabilization radius larger than $M_n$ is bounded, i.e.,
\begin{equation}\label{inequality_nodes_nearest_neighbor_large_stabilization}
\#\mc J_n^M \le k 2^d n / (\kappa_d M^d).
\end{equation}
This can be seen by going through a configuration from $\PP_n$ node by node and assigning the labels essential and inessential to some of them. Each considered node $X\in\PP_n$ with stabilization radius larger than $M_n$ that has not been labeled yet, is labeled as essential and each of its $k-1$ closest neighbors is labeled as inessential if it has not been labeled as essential before. After the procedure, all essential nodes cannot have any other essential points within distance $M_n$. Consequently, balls with radius $M_n/2$ around the essential nodes cannot intersect. The bound in \eqref{inequality_nodes_nearest_neighbor_large_stabilization} is derived by bounding the number of these balls in $[0,1]^d$ through the volume each occupies and multiplying with $k$ to adjust for the inessential points.

An issue that can arise when it comes to the sprinkling requirements are relatively close nodes in $\mc J_n^M$, due to potentially not disjoint sets in the sprinkling event. To make sure that such scenarios cannot occur, we aim to thin out the set of these bad vertices. We say a node $X\in\PP_n$ is \emph{distinguished} if $X$ is the smallest node in the lexicographic order of $(\PP_n\cap B_{n^{-1/d}}(X))-X$. Then, we define
$$
\widetilde{\mc J}_n^M := \widetilde{\mc J}_n^M(\PP_n) := \{X\in\mc J_n^M\colon X\text{ is distinguished}\}
$$
as a subset of $\mc J_n^M$ that only keeps distinguished nodes. This guarantees that the distance between two nodes in $\widetilde{\mc J}_n^M$ is at least $n^{-1/d}$ and therefore balls with radius $n^{-1/d}/2$ centered in each node in $\widetilde{\mc J}_n^M$ are disjoint. With this in mind, we can define the sprinkling event by setting
$$F_n^{M,(2)} := \big\{\PP_n^{+,M}\big([0,1]^d \setminus \cup_{X\in\widetilde{\mc J}_n^M} B_{n^{-1/d}/2}(X)\big) = 0\big\}$$
and
$$F_n^{M,(3)} := \bigcap_{X\in\widetilde{\mc J}_n^M} \{\PP_n^{+,M}(B_{n^{-1/d}/2}(X)) = k\}.$$
Further, if we assume that $M$ is large, it follows that each node $X\in\mc J_n^M$ can only have $k-1$ other nodes in $B_{k n^{-1/d}}(X)$. Otherwise, $X$ would have a stabilization radius bounded by $k n^{-1/d}$. Thus, one of the nodes in $B_{k n^{-1/d}}(X)$ has to be distinguished. Subsequently, after adding $k$ points to $B_{n^{-1/d}/2}(X)$ for each $X\in\widetilde{\mc J}_n^M$, the stabilization radius for each $X\in\mc J_n^M$ is bounded by $(k+1)n^{-1/d}$ and the same bound holds for the stabilization radii of the additionally inserted points. Hence, \eqref{SPR1} is fulfilled and the bound from \eqref{inequality_nodes_nearest_neighbor_large_stabilization} with the definition of $\widetilde{\mc J}_n^M$ implies \eqref{SPR2} with $I_n^M(\PP_n) := \#\widetilde{\mc J}_n^M$, $V:= 2^{-d}$ and $m:=k$.

Finally, to verify \eqref{SPR3}, we see that the $k$ nodes put in $B_{n^{-1/d}/2}(X)$ for every $X\in\widetilde{\mc J}_n^M$, each come with an additional score that is bounded by $k$ after the rescaling with $n^{1/d}$. All scores of vertices that already existed can only decrease when inserting the new nodes and the same holds for the cut-off score. This means we arrive at
$$H_n^{M',M} (\PP_n^M) \leq H_n^{M',M}(\PP_n) + \frac1{n} k^2 \#\mc J_n^M \le  H_n^{M',M}(\PP_n) + \frac{\log M}{M^d}$$
under $F_n^M$, for large $M$ and $M'>M$, confirming \eqref{SPR3}.

\subsubsection{Power-weighted spherical contact distances}
A basic characteristic of a point pattern is the distribution of the spherical contact distances \cite[Section 4.2]{illian}. Loosely speaking, it describes the distance to the nearest point of the given point pattern measured from a space point that is selected at random. A basic approach to estimate this quantity is the point-count method \cite[Section 4.2]{illian}. Here, the window is discretized, and then for each subcube, the distance of its center to the closest point is recorded. A natural way to formulate an estimator that is independent of the discretization, is to replace the discretization with an integral. Following this setup, in the present example, we describe the large deviation behavior of estimators of the $\a$th moment of the spherical contact distances for $\a > d$. For this, we aim to use the integral form representation \textbf{B}. We define the score function by
$$\xi(x,\varphi) :=  \inf_{y\in\varphi} \|y - x\|^\alpha,$$
where, as in the previous section, we replace the Euclidean norm with the toroidal distance when applying the score function to a configuration on a torus.

\begin{figure}[H]
\centering
\begin{tikzpicture}[scale = 5.85]
\input{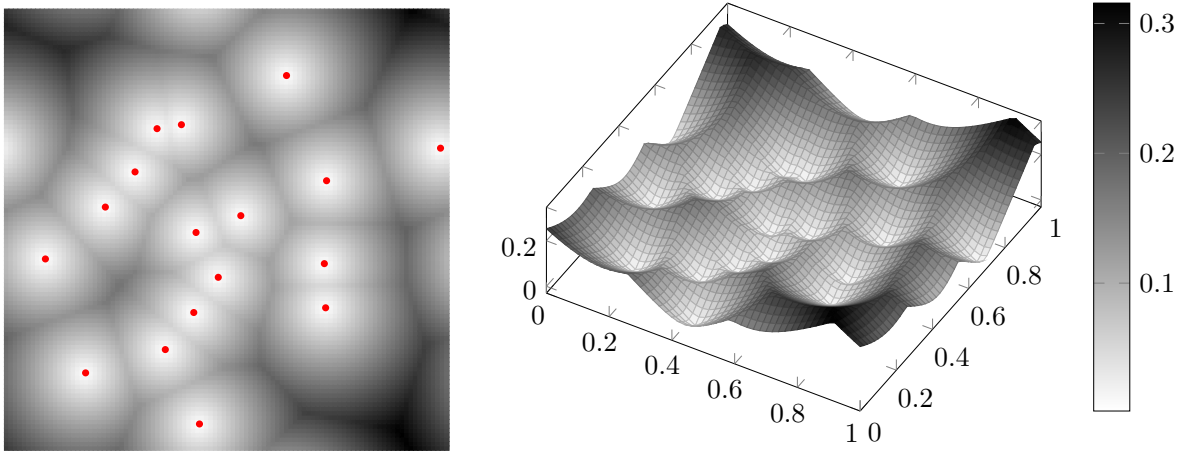}
\end{tikzpicture}
\hspace{.3cm}
\begin{tikzpicture}[scale = .95]
\begin{axis}[
  colorbar,
  colormap name = whiteblack,
  view={30}{70},
]
\addplot3[surf, mesh/rows=51] table [x index=0, y index=1, z index=2] {Tikz_pictures/contact_distances.txt};
\end{axis}
\end{tikzpicture}
\caption{Simulation of spherical contact distances based on a Poisson point process on a two-dimensional torus. The lighter the shade, the smaller the distance of a space point to its closest node in the configuration.}\label{fig:contact_distances}
\end{figure}

We can set the stabilization radius to be
$$\mc R(x,\varphi) := \inf\{r > 0 \colon \varphi(B_r(x)) \ge 1\}.$$
With this stabilization radius and the choice $g(M):=M^\alpha$, \eqref{STA} is satisfied. Also, \eqref{INC} holds with the choice $\xi^r(x,\varphi) := \xi(x,\varphi\cap B_r(x))\wedge r$ for $r>0$, $x\in\R^d$ and $\varphi\in\NNNN$.

In order to construct the sprinkling event, we divide $[0,1]^d$ into a grid of cubes of side length $n^{-1/d} M/\log M$, denote this collection by $\mc Q_n^M$ and call a box $Q\in\mc Q_n^M$ bad if it does not contain any Poisson points, i.e., if $Q\cap\PP_n = \emptyset$. Let
$$\mc J_n^M := \mc J_n^M(\PP_n) := \{Q\in\mc Q_n^M\colon \PP_n\cap Q=\emptyset\}$$
be the set of bad boxes. If a box is bad, all points in a cube of volume $n^{-1/d}$ in the center of the bad box must have a distance to the closest node of at least $n^{-1/d} M/(\log M)^2$ for sufficiently large $M$. Thus, a bad subcube contributes with a value of at least $M^\alpha /(\log M)^{2\alpha}$ to the total functional after resolving the rescaling with factor $n^{1/d}$. Thus, under the event $\{H_n^{M',M} < a\}$, for $M'>M$, such bad boxes can only occur a limited number of times. More precisely, due to our choice of $g$, we find that
\begin{equation}\label{equation_bad_boxes_contact_distance1}
\#\mc J_n^M \le n a /(M^\alpha /(\log M)^{2\alpha}) = n a \frac{(\log M)^{2\alpha}}{M^\alpha}.
\end{equation}
Now, to define the sprinkling event, we introduce an additional sub grid. First, without explicitly stating it, in the following, we will assume that $M$ is sufficiently large for some properties to hold and that we can manage the assignment of the subcubes without having to deal with fractions of subcubes. A negligible adjustment of the side length of the boxes would assure the latter. Divide $Q\in\mc Q_n^M$ into subcubes of side length $n^{-1/d} \log M$ and call this collection $\mc W_n^M(Q)$. With the observation in \eqref{equation_bad_boxes_contact_distance1}, the number of subcubes in bad boxes is bounded by
\begin{equation}\label{equation_bad_boxes_contact_distance2}
\#\{W\in \mc W_n^M(Q) \colon Q\in \mc J_n^M\} \le n a \frac{(\log M)^{2\alpha}}{M^\alpha} \frac{(n^{-1/d} M/\log M)^d}{(n^{-1/d} \log M)^d} = n a \frac{(\log M)^{2(\alpha-d)}}{M^{\alpha-d}}.
\end{equation}
For $Q\in\mc J_n^M$ and $W\in\mc W_n^M(Q)$, let $B_W\subseteq W$ be the ball with radius $n^{-1/d}$ that is located around the center of $W$. Now, we can define the sprinkling event by inserting a node in each subcube of every bad box. Thus, accordingly to \eqref{equation_sprinkling_event1} and \eqref{equation_sprinkling_event2}, we get the events
$$F_n^{M,(2)} := \big\{\PP_n^{+,M}\big([0,1]^d \setminus (\cup_{W\in\{W'\in \mc W_n^M(Q) \colon Q\in \mc J_n^M\}} B_W)\big) = 0\big\}$$
and
$$F_n^{M,(3)} := \bigcap_{W\in\{W'\in \mc W_n^M(Q) \colon Q\in \mc J_n^M\}} \{\PP_n^{+,M}(B_W) = 1\}.$$
If all bad boxes contain at least one vertex, the stabilization radius of any space point can be at most of order $n^{-1/d} M/\log M$, and is therefore, less than or equal to $M/n^{1/d}$ for large enough $M$, verifying \eqref{SPR1}. Further, since we assumed $\alpha>d$, \eqref{equation_bad_boxes_contact_distance2} confirms \eqref{SPR2} with $I_n^M(\PP_n) := \#\{W\in \mc W_n^M(Q) \colon Q\in \mc J_n^M\}$, $V:= 1$ and $m:=1$. For \eqref{SPR3}, we point out that inserting an additional node cannot increase the contact distance of any point. Additionally, every space point in a good box has a contact distance of order $n^{-1/d} M /\log M$ and, thus, cannot be affected by the cut-off of the score in the functional and thus, also with respect to the cut-off functional the contact distance of a space point in a good box after the sprinkling can only decrease. This observation yields that only the added points in bad boxes have to be considered to bound the increase of the cut-off functional under the sprinkling event. But under $F_n^{M,(3)}$, the distance to the closest node of every space point in a bad box is of order $\log M$ and thus, bounded by $(\log M)^2$ for large $M$. Hence, we get that under $F_n^M\cap \{H_n^{M',M} < a\}$
$$
H_n^M (\PP_n^M) = H_n^{M',M}(\PP_n^M) \le H_n^{M',M}(\PP_n) + \frac{(\log M)^{2\alpha}}{n} I_n^M(\PP_n) \le H_n^M(\PP_n) + a \frac{(\log M)^{4\alpha-2d}}{M^{\alpha-d}},
$$
also verifying \eqref{SPR3}.

\subsection{Functionals for the sparse random geometric graph}\label{section_applications_sparse}
\subsubsection{Subgraph counts}
Let $G_0 := (V,E)$, where $V$ represents a set of vertices and $E$ a set of edges, be an arbitrary fixed finite connected graph. With this, we define
$$\xi(\varphi) := \#\big\{(\varphi', E') \colon \varphi'\subseteq\varphi,\; E'\subseteq\{\{x,y\}\subseteq \varphi'\colon \|x-y\|\le 1\},\; (\varphi', E') \cong G_0\big\}$$
for a configuration $\varphi\in\NNN$, to count the occurrence of the graph $G_0$ in the geometric graph with connectivity radius $1$ on $\varphi$. For $n\in\N$ and $\varphi\in\NNN^{(1)}$ with $\dist(\varphi,\partial[0,1]^d) > r_n$, we define $\xi_n(\varphi)$ similar to $\xi(r_n^{-1}\varphi)$ but replace the Euclidean distance with the toroidal metric of $[0,r_n^{-1}]^d/\sim$. These functionals fulfill all requirements stated in Theorem \ref{theorem_main_sparse}. If used as a score function, as displayed in \eqref{equation_functional_RGG}, it represents occurrences of $G_0$ in a random geometric graph with connectivity radius $r_n$ in a sparse regime.

Additionally, sometimes it is possible to simplify the rate function further. More precisely, assume that we count the occurrences of a $k_0$-clique. Then, with Mecke's formula, it can be computed that
\begin{equation}\label{equation_example_subgraph_count_expectation1}
\E[H_n^\sp] \overset{n\to\infty}{\longrightarrow} \frac{v_{d,k}(G_0)}{k_0!} =: \mu_{d,k_0}.
\end{equation}
The right-hand side is given  by
\begin{align}\label{equation_example_subgraph_count_expectation2}
\begin{split}
v_{d,k_0}(G_0) &:= \int_{\R^{d(k_0-1)}} \prod_{\{i, j\}\in \{1,\dots,k_0\}} \one\big\{\|x_i-x_j\| \le 1 \big\} \d(x_2,\dots,x_{k_0}),
\end{split}
\end{align}
where $x_1:=0$ and $v_{d,k_0}(G_0):=1$ if $k_0=1$. Intuitively, \eqref{equation_example_subgraph_count_expectation2} represents the volume of all possible locations to place $k_0-1$ points around a fixed point such that the generated geometric graph with connectivity radius $1$ is isomorphic to $G_0$. Now, from our proof for the sparse regime, it follows that we can also write $\xi$ directly as an indicator that triggers for complete connected components of size $k_0$. Then, \cite[Remark 3.6]{hirschowada} implies that
\begin{align*}
-\inf_{\rho\colon T^\sp(\rho) \le a} h_{k_0}^\sp(\rho \mid \tau_{k_0}^\sp) &= -\inf_{x \le a} x\log(x/\mu_{d,k_0}) - x + \mu_{d,k_0} \\
&= \begin{cases} -a\log(a/\mu_{d,k_0}) + a - \mu_{d,k_0} &\text{if } a<\mu_{d,k_0} \\ 0 &\text{otherwise}\end{cases},
\end{align*}
and an analogous simplification could be achieved but would require substantial additional computations and is therefore omitted.

\subsubsection{Betti numbers and persistent Betti numbers}
Simply expressed, Betti numbers count holes of a certain dimension in simplicial complexes. \cite[Section 4.1]{hirschowada} gives a short overview of literature dealing with the basic concepts behind Betti numbers and the more general persistent Betti numbers. They can be built upon the \v{C}ech complex. For a set $\varphi\in\NNN$ and $r\ge0$, the \v{C}ech complex is defined by
$$\check{C}_r(\varphi) := \big\{\psi\subseteq\varphi\colon \cap_{x\in\psi} B_r(x) \neq \emptyset\big\}.$$
Now, we can define the $k$th persistent Betti number for $0\le s\le t\le \infty$ by
$$\beta_k(\varphi,s,t) = \dim \frac{Z_k(\check{C}_s(\varphi))}{Z_k(\check{C}_s(\varphi))\cap \widetilde B_k(\check{C}_t(\varphi))},$$
where $\varphi\in\NNN$, $Z_k$ is the $k$th cycle group of the \v{C}ech complex and $\widetilde B_k$ represents the $k$th boundary group. For configurations close to $\partial[0,1]^d$, we define $\xi_n$ similar to $\xi(r_n^{-1}\,\cdot\,)$, using balls with respect to the torus $[0,r_n^{-1}]^d/\sim$ to set up the \v{C}ech complex. The requirements for Theorem \ref{theorem_main_sparse} are satisfied and we recover the ordinary Betti numbers by setting $s=t$. As in the case of subgraph counts, also here, a simplification of the rate function according to \cite[Remark 3.6]{hirschowada} is achievable. However, to keep this section at a reasonable size, we omit the explicit computations.

\subsubsection{Edge lengths}
For a point set $\varphi\in\NNN$, we define
$$\xi(\varphi) := \sum_{x,y\in\varphi} \|x-y\| \one\{\|x-y\| \le 1\}$$
and $\xi_n$, for $n\in\N$, is defined analogously to the subgraph counts or Betti numbers examples, using the toroidal metric of $[0,r_n^{-1}]^d/\sim$ instead of the Euclidean distance.
Then, all requirements of Theorem \ref{theorem_main_sparse} are satisfied. Note that here, $k_0=2$ and thus, as the proof of Theorem \ref{theorem_main_sparse} shows, only isolated edges will be relevant for the lower large deviations.
\section{Proof of Theorem \ref{theorem_main_thermodynamic} (critical)}\label{section_proof_thermodynamic}

For bounded and local score functions, \cite{georgii} provides a large deviation principle for associated functionals. We recall that our strategy is to use a coupling consisting of a thinned Poisson point process and another independent Poisson process. We let the thinning fully replicate $\PP_n$ while using the independent Poisson point process to sprinkle in additional points following a specific pattern to guarantee locality and boundedness of the score function such that the general large deviations theory becomes invokable.

First, we let $\widetilde\PP_n$ be a Poisson point process with intensity $1$ on the torus $[0,n^{1/d}]^d / \sim$. Note that $n^{1/d}\PP_n$ and $\widetilde\PP_n$ have the same distribution. Now, we can replicate the proof of \cite[Theorem 1.1]{hirsch} to get Theorem \ref{theorem_main_thermodynamic} a), the upper bound for the lower large deviations.

\bep[Proof of Theorem \ref{theorem_main_thermodynamic} a)]
We recall $\xi^r$ from \eqref{INC}. Further, for the next steps, we assume that we are in representation \textbf{A} and indicate that the other case works analogously. The functional $\xi^r$ is bounded and local, and thus, we can use \cite[Theorem 3.1]{georgii} (or \cite[Corollary 3.2]{georgii} in the case of representation \textbf{B}) to get that
\begin{align*}
\limsup_{n\to\infty} \frac1{n} \log\P(H_n^\th \le a) &\le \limsup_{n\to\infty} \frac1{n} \log\P\bigg(\frac1{n}\sum_{X\in\PP_n} \xi^r(n^{1/d}X, n^{1/d}\PP_n) \le a\bigg) \\
&= \limsup_{n\to\infty} \frac1{n} \log\P\bigg(\frac1{n}\sum_{X\in\widetilde\PP_n} \xi^r(X, \widetilde\PP_n) \le a\bigg) \le -\inf_{\Q\colon\Q^o[\xi^r] \le a} h^\th(\Q).
\end{align*}
By \eqref{INC}, $\xi^r(x,\varphi)$ increases, as $r$ grows, towards $\xi(x,\varphi)$ for each $x\in\varphi\in\NNNN$. Proceeding, using monotone convergence, as in the proof of \cite[Theorem 1.1]{hirsch}, it follows that
$$-\limsup_{r\to\infty} \inf_{\Q\colon\Q^o[\xi^r] \le a} h^\th(\Q) \le -\inf_{\Q\colon\Q^o[\xi] \le a} h^\th(\Q),$$
which concludes the upper bound.
\enp

In order to prove the lower bound, it is necessary to examine the event $F_n^M$ from \eqref{SPR1} in detail. For this, we denote the number of Poisson points of $\PP_n$ by $N_n := \PP_n([0,1]^d)$. The next lemma gives a lower bound for the probability of the sprinkling event.

\bel[Sprinkling regularizes with high probability]\label{lemma_sprinkling} For $n\geq M\geq 1$ sufficiently large, we get that almost surely
$$\P(F_n^M\mid\PP_n) \geq (1-M^{-1})^{N_n} e^{-n/M} \big( \tfrac{(V/M)^m}{m!} e^{-V/M}\big)^{I_n^M(\PP_n)}.$$
\enl

\bep[Proof of Lemma \ref{lemma_sprinkling}]
Looking at the probabilities of each single event of $F_n^M$ gives
\begin{align*}
&\P(\PP_n^{-,M}= \PP_n \mid \PP_n) = (1-M^{-1})^{N_n}, \\
&\P\big(\PP_n^{+,M}\big([0,1]^d \setminus (\cup_{i = 1}^{I_n^M(\PP_n)} B_{n,i}^M(\PP_n)) = 0 \mid \PP_n\big) \geq e^{-n/M}
\intertext{and}
&\P\bigg(\bigcap_{i = 1}^{I_n^M(\PP_n)} \{\PP_n^{+,M}(B_{n,i}^M(\PP_n)) = m\} \biggm\vert \PP_n\bigg) = \big(\tfrac{(V/M)^m}{m!} e^{-V/M}\big)^{I_n^M(\PP_n)}
\end{align*}
almost surely, where we used that the survival probability of the thinning is $1-M^{-1}$ and the intensity of $\PP_n^{+,M}$ was assumed to be $n/M$. Using independence between all three events conditioned on $\PP_n$ yields the desired statement.
\enp

Now, we conclude the proof of Theorem \ref{theorem_main_thermodynamic}.

\bep[Proof of Theorem \ref{theorem_main_thermodynamic} b)]
In the following, assume that $M>0$ is large and $M'>M$. Because of $F_n^M\subseteq E_n^M$, which was assumed in \eqref{SPR1}, it follows that under the event $F_n^M$ the radius of stabilization with respect to $\PP_n^M$ of each node in $\PP_n^M$ or space point in $[0,1]^d$, depending on whether we consider a functional given in representation \textbf{A} or representation \textbf{B}, is at most $M n^{-1/d}$. Hence, we can invoke \eqref{STA} from which follows that under $F_n^M$ we can replace $H_n(\PP_n^M)$ by $H_n^M(\PP_n^M)$, and get
$$\P(H_n(\PP_n)<a) = \P(H_n(\PP_n^M)<a) \geq \P(\{H_n^M(\PP_n^M) < a\} \cap F_n^M\cap\{H_n^{M',M}(\PP_n) < a\}).$$
Due to \eqref{SPR3} it holds that under $F_n^M\cap\{H_n^{M',M}(\PP_n) < a\}$ almost surely
\begin{equation}\label{inequality_added_points}
H_n^M (\PP_n^M) \leq H_n^{M',M}(\PP_n) + c_M^{(2)},
\end{equation}
as $M\rightarrow\infty$ and thus, since $\{H_n^{M',M}(\PP_n) < a - c_M^{(2)}\} \subseteq \{H_n^{M',M}(\PP_n) < a\}$, it follows that,
$$\P(\{H_n^M(\PP_n^M) < a\} \cap F_n^M \cap\{H_n^{M',M}(\PP_n) < a\}) \geq \P(\{H_n^{M',M}(\PP_n) + c_M^{(2)} < a\} \cap F_n^M).$$
By conditioning on $\PP_n$ and applying Lemma \ref{lemma_sprinkling} for sufficiently large $n$, we arrive at
\begin{align*}
&\E\big[\one\{H_n^{M',M}(\PP_n) < a - c_M^{(2)}\} \P(F_n^M\mid\PP_n)\big] \\
&\ge \E\big[\one\{H_n^{M',M}(\PP_n) < a - c_M^{(2)}\} (1-M^{-1})^{N_n}  \big(\tfrac{(V/M)^m}{m!} e^{-V/M}\big)^{I_n^M(\PP_n)}\big] e^{-n/M} \\
&= \E\Big[\one\{H_n^{M',M}(\PP_n) < a - c_M^{(2)}\} \exp\Big(N_n\log(1-M^{-1}) + I_n^M(\PP_n)\log\big(\tfrac{(V/M)^m}{m!} e^{-V/M}\big)\Big)\Big] e^{-n/M}.
\end{align*}
Moreover, invoking \eqref{SPR2} and introducing a bound for $N_n$ yields for any $c>0$,
\begin{align*}
&\E\Big[\one\{H_n^{M',M}(\PP_n) < a - c_M^{(2)}\} \exp\Big(N_n\log(1-M^{-1}) +  I_n^M(\PP_n)\log\big(\tfrac{(V/M)^m}{m!} e^{-V/M}\big)\Big)\Big] \\
&\ge \P\big(H_n^{M',M}(\PP_n) < a - c_M^{(2)}, N_n < cn\big) \exp\big(c n\log(1-M^{-1}) + c_M^{(1)} n \log(\tfrac{(V/M)^m}{m!} e^{-V/M})\big).
\end{align*}
To convince ourselves that the exponential factors are not relevant, we recall that $c_M^{(1)}/\log M \rightarrow 0$ as $M\rightarrow\infty$ was assumed, which yields
\begin{align*}
&\frac1{n}\log\Big(\exp\big(cn \log(1-M^{-1}) + c_M^{(1)} n \log(\tfrac{(V/M)^m}{m!} e^{-V/M}) - n/M\big)\Big) \\
&= c \log(1-M^{-1}) + c_M^{(1)} \big(\log(\tfrac{(V/M)^m}{m!}) - V/M\big) - M^{-1} \overset{M\to\infty}{\longrightarrow} 0.
\end{align*}
Now, for the other factor,
\begin{align*}
&\P\big(H_n^{M',M}(\PP_n) < a - c_M^{(2)}, N_n < c n\big) \ge \P\big(H_n^{M',M}(\PP_n) < a - c_M^{(2)}\big) - \P\big(N_n \ge c n\big),
\end{align*}
where for large $c$, \cite[Lemma 1.2]{poisson_conc} can be used to show that the second term does not affect the large deviations.

For the next computations, we assume that $H_n^\th$ has representation \textbf{A}. The other case works analogously.  We define $\xi^{M',M}(x,\varphi) := \xi(x, \varphi \cap B_{M'}(x))\wedge g(M)$ for $x\in\R^d$ and $\varphi\in\NNNN$, and point out that $\xi^{M',M}$ can be locally determined and is bounded by $g(M)$. Besides that, recall that $\widetilde{\PP}_n$ is equal in distribution to $n^{1/d}\PP_n$. Then, applying \cite[Theorem 3.1]{georgii} (or \cite[Corollary 3.2]{georgii} in the case of representation \textbf{B}), we can proceed as in the proof of the upper bound, and we arrive at
\begin{align*}
&\liminf_{n\to\infty} \frac1{n} \log \P\big(H_n^{M',M}(\PP_n) < a - c_M^{(2)}\big) \\
&= \liminf_{n\to\infty} \frac1{n} \log\P\Big(\frac1{n} \sum_{X\in\PP_n} \xi(n^{1/d}X, n^{1/d}\PP_n \cap B_{M'}(n^{1/d} X))\wedge g(M) \le a - c_M^{(2)}\Big) \\
&= \liminf_{n\to\infty} \frac1{n} \log\P\Big(\frac1{n} \sum_{X\in\widetilde{\PP}_n} \xi^{M',M}(X,\widetilde\PP_n) \le a - c_M^{(2)}\Big) \ge -\inf_{\Q\colon\Q^o[\xi^{M',M}]<a - c_M^{(2)}} h^\th(\Q).
\end{align*}
Finally, we assert that
$$\liminf_{M\rightarrow\infty} \liminf_{M'\rightarrow\infty} \Big(-\inf_{\Q\colon\Q^o[\xi^{M',M}]<a - c_M^{(2)}} h^\th(\Q)\Big) \ge -\inf_{\Q\colon\Q^o[\xi]<a} h^\th(\Q),$$
which yields the desired result.

To prove this assertion, let $\Q$ be an arbitrary point process that satisfies $\E_{\Q^o}[\xi(0,\cdot\, )]<a$. This lets us find some $\delta>0$ such that $\E_{\Q^o}[\xi(0,\cdot\, )]<a-\delta$. Next, for any $M>0$, it also holds that $\E_{\Q^o}[\xi(0,\cdot\, ) \wedge g(M)]<a-\delta$ due to monotonicity. Further, dominated convergence yields that $\lim_{M'\to\infty} \E_{\Q^o}[\xi(0,\cdot\, \cap B_{M'}(0)) \wedge g(M)] = \E_{\Q^o}[\xi(0,\cdot\, ) \wedge g(M)]$ from which we deduce the existence of $M_0(\delta,M)$ such that for all $M'>M_0(\delta,M)$
$$\E_{\Q^o}[\xi(0,\cdot\, \cap B_{M'}(0)) \wedge g(M)] < a - \delta/2.$$
In particular, from $c_M^{(2)} \rightarrow 0$ as $M\rightarrow\infty$ we get that for some $M_0(M)>0$ and all $M'>M_0(M)$
$$\E_{\Q^o}[\xi(0,\cdot\, \cap B_{M'}(0)) \wedge g(M)] < a - c_M^{(2)}$$
if $M$ is large enough. Therefore,
$$\{\Q\colon\Q^o[\xi]<a\} \subseteq \{\Q\colon\Q^o[\xi^{M',M}]<a - c_M^{(2)} \text{ for all }M'>M_0(M) \text{ and }M\text{ large}\}$$
which implies that
$$\limsup_{M\rightarrow\infty} \limsup_{M'\rightarrow\infty} \Big(\inf_{\Q\colon\Q^o[\xi^{M',M}]<a - c_M^{(2)}} h^\th(\Q)\Big) \le \inf_{\Q\colon\Q^o[\xi]<a} h^\th(\Q).$$
\enp

\section{Proof of Theorem \ref{theorem_main_sparse} (sparse)}\label{section_proof_sparse}

For the sparse case, we would like to apply the large deviation principle for empirical measures counting potentially connected components of a fixed size of a random geometric graph from \cite[Theorem 2.1]{hirschowada}. Using sprinkling, we would ideally like to create a coupled Poisson point process that, when serving as nodes for a geometric graph, only contains fixed-sized components. A simple replication of the procedure in the critical case for the sparse case is not possible as we will desire for the thinning to keep most of the points, which will be with very high probability an amount of order $n$, thus, resulting in costs for the thinning of magnitude $e^{-c n}$ for some $c>0$. But the speed for the sparse regime satisfies
$$\frac{n}{\rho_{n,k_0}^\sp} = \frac{n}{n^{k_0} r_n^{d(k_0-1)}} = (n r_n^d)^{-(k_0-1)} \rightarrow \infty,$$
if $k_0>1$. Instead, as in \cite{hirschowada}, we will divide $[0,1]^d$ into a grid and resample an entire box of the grid if we deem the configuration in it as not feasible and additionally bound the inevitable error in the functional that this process creates. This then results in a coupled Poisson process as a foundation for a geometric graph for which all significant connected components are of a fixed size, and therefore, we can invoke the large deviation principle from \cite[Theorem 2.1]{hirschowada}. 

Now, to give more details after this overview, as announced, we start by dividing $[0,1]^d$ into a grid of cubes with side length $(\rho_{n,k_0}^\sp)^{-1/d}$ each, where to keep the notation simpler, we assume that $\rho_{n,k_0}^\sp$ is a natural number and denote this collection by $\mathcal{Q}_n$. We define $\PP'_n$ as an Poisson point process on $[0,1]^d$ with intensity $n$ independent of $\PP_n$. Further, for all cubes $Q\in\mathcal{Q}_n$, let $X_{Q, \eps}$ be Bernoulli random variables with parameter $\eps\in(0,1)$, independent of each other and all introduced Poisson random measures. Using this, we define
$$\PP_n^Q := \begin{cases}
	Q\cap\PP'_n &\text{if } X_{Q, \eps} = 1 \\
	Q\cap\PP_n &\text{if } X_{Q, \eps} = 0
	\end{cases},$$
which yields a Poisson point process on $Q$ with intensity $n$ for each $n\in\N$. Consequently, $\PP''_n := \cup_{Q\in\mathcal{Q}_n} \PP_n^Q$ is a Poisson point process on $[0,1]^d$ with intensity $n$. The idea is to use the Bernoulli random variables to control $\PP''_n$ in such a way that we resample $\PP_n$ using $\PP'_n$ in each box that has a node with $k_0$ relatively close other vertices while keeping $\PP_n$ in all other boxes.
To achieve this, let
$$\mathcal{J}_n := \mathcal{J}_n(\PP_n) := \big\{Q \in \mathcal{Q}_n\colon \max_{X\in Q\cap\PP_n} \PP_n(B_{2^d k_0 r_n}(X)) \ge k_0+1\big\}$$
be the boxes that contain a vertex with at least $k_0$ other vertices within distance $2^d k_0 r_n$ and that we would therefore like to resample. To further ease notation, we also denote the number of bad boxes by
$$J_n := \#\mathcal{J}_n(\PP_n)$$
and we point out that we can consider $k_0$ as fixed from now on, which lets us write
$$\rho_n^\sp := \rho_{n,k_0}^\sp.$$
We first make sure that these bad boxes do not occur too many times with a probability that is too high.

\bel[Bad boxes are exponentially negligible]\label{lemma_probability_bad_boxes_sparse}
Let $\delta>0$. Assume that $n r_n^d\rightarrow 0$ and $\rho_n^\sp \rightarrow\infty$. Then,
$$\limsup_{n\to\infty} \frac1{\rho_n^\sp} \log \P(J_n \ge \delta\rho_n^\sp) = -\infty.$$
\enl

Next, we determine what happens within a box that was resampled and ignore effects of adjacent boxes for now. Preferably we would like the sprinkled process not to create any new components consisting of $k_0$ or more vertices within a resampled cube. The next lemma states that for each $n\in\N$, conditioned on $\PP_n$, the probability of not having $k_0$ close points within a resampled box $Q\in\mathcal{J}_n$ is bounded from below.

\bel[With positive probability, a resampled box does not contain $k_0$ close nodes]\label{lemma_prob_new_edges_sparse}
Assume that $n r_n^d\rightarrow 0$ and $\rho_n^\sp\rightarrow\infty$. Then, for any $M > 2 \kappa_d^{k_0-1} 2^{k_0(d^2+1)} k_0^{k_0}$ it holds that
$$\P\bigg(\bigcap_{Q\in\mathcal{J}_n} \Big\{\max_{X\in Q\cap\PP'_n}  \PP'_n\big(B_{2^d k_0 r_n}(X)\cap Q\big) \le k_0-1\Big\} \biggm\vert \PP_n\bigg) \ge \alpha_M^{J_n},$$
where $\alpha_M := 2^{-M-1} $.
\enl

One issue that we cannot prevent is that there can be large connected components between two adjacent boxes when at least one of them is resampled. But we can show that the number of these components will, with high enough probability, not be significant. More precisely, the next lemma will control the number of large components that can occur between boxes when resampling. To ease notation, for every $Q\in\mathcal{Q}_n$, we let
$$\partial_n Q := \{x\in Q\colon \dist(\{x\},\partial Q) \le 2^d k_0 r_n\}$$
denote the set of all points in $Q$ within distance $2^d k_0 r_n$ of the boundary of $Q$. The factor $2^d k_0$ appears here to be able to deal with boxes that share a face, which allows for large connected components to exist that span over multiple boxes. We also let
$$\CCC_{n,k_0} := \{X\in\PP''_n\colon s_n(\{X\}\cup\varphi,\PP''_n)=1 \text{ for some }\varphi\subseteq\PP''_n \text{ with } \#(\{X\}\cup\varphi)\in{\{k_0,\dots, 2^d k_0\}}\}$$
be the vertices in $\PP''_n$ that are part of a connected component of size between $k_0$ and $2^d k_0$, where we recall the definition of $s_n$ from \eqref{equation_indicator_torus}.

\bel[The number of large connected components between boxes is negligible]\label{lemma_prob_resampling_between_boxes_sparse}
Let $\delta>0$. Assume that $n r_n^d\rightarrow 0$ and $\rho_n^\sp\rightarrow\infty$. Then,
$$\limsup_{n\to\infty} \frac1{\rho_n^\sp} \log \P\big(\CCC_{n,k_0}(\cup_{Q\in \mathcal{Q}_n} \partial_n Q) \ge \delta\rho_n^\sp\big) = -\infty.$$
\enl

With these lemmas and preliminaries, we can prove the lower large deviations in the sparse regime.

\bep[Proof of Theorem \ref{theorem_main_sparse}]
We point out that \cite{hirschowada} worked with the Euclidean distance on $[0,1]^d$ instead of the toroidal metric. For this reason, we need some additional notation to deal with this subtle difference. Also, recall that $H_n^\sp = \frac{1}{\rho_n^\sp} \sum_{\varphi\subseteq \PP_n} \xi_n(\varphi) s_n(\varphi,\PP_n)$, where $s_n$ checks for connected components with respect to the toroidal metric. We define the restriction to components of size $k_0$ by
$$H_{n,k_0}^\sp := H_{n,k_0}^\sp(\PP_n) :=\frac{1}{\rho_n^\sp} \sum_{\varphi\subseteq \PP_n, \#\varphi = k_0} \xi_n(\varphi) s_n(\varphi,\PP_n).$$
To also incorporate the Euclidean metric, we define $\ms{GG}'_n(\varphi)$ as the geometric graph on $\varphi\in\NNN$ with connectivity radius $r_n$ and with respect to the Euclidean distance. With this, for $\varphi\subseteq\psi\in\NNN$, we set
$$s'_n(\varphi,\psi) := \one\{\varphi\text{ is a connected component of }\ms{GG}'_n(\psi)\}$$
to be the counterpart of $s_n$ in terms of the Euclidean distance. Further, along the lines of \cite[Theorem 3.3]{hirschowada}, we define
$$\widetilde{H}_{n,k_0}^\sp := \widetilde{H}_{n,k_0}^\sp(\PP_n) := \frac1{\rho_n^\sp} \sum_{\varphi\subseteq\PP_n, \#\varphi=k_0} \xi(r_n^{-1}\varphi) t_n(\varphi,\PP_n),$$
where
$$t_n(\varphi,\PP_n) := \one\{\|y-z\| \ge r_n \text{ for all }y\in \varphi \text{ and } z\in \PP_n\setminus \varphi\} \one\{\diam(\varphi)\le k_0 r_n\}$$
is the indicator assuring that $\varphi$ is isolated and locally concentrated within $\PP_n$.

Our goal is to apply \cite[Theorem 3.3]{hirschowada} to $\widetilde{H}_{n,k_0}^\sp$. One main step for the upper bound of this proof will be to show that the error between $H_{n,k_0}^\sp$ and $\widetilde{H}_{n,k_0}^\sp$ that occurs close to the boundary is negligible.
Thus, we define
$$H_{n,k_0}^{\text{err},1}(\PP_n) := \frac{1}{\rho_n^\sp} \sum_{\substack{\varphi\subseteq\PP_n\colon \#\varphi=k_0, \\ \dist(\varphi,\partial[0,1]^d) \le r_n}} \xi(r_n^{-1}\varphi) s'_n(\varphi,\PP_n)$$
and compute
\begin{align}\label{inequality_fixed_size_error_sparse}
\begin{split}
H_n^\sp &\ge H_{n,k_0}^\sp \ge \frac{1}{\rho_n^\sp} \sum_{\substack{\varphi\subseteq\PP_n\colon \#\varphi=k_0, \\ \dist(\varphi,\partial[0,1]^d) > r_n}} \xi(r_n^{-1}\varphi) s'_n(\varphi,\PP_n) \\
&= \frac{1}{\rho_n^\sp} \sum_{\varphi\subseteq\PP_n, \#\varphi=k_0} \xi(r_n^{-1}\varphi) t_n(\varphi,\PP_n) - \frac{1}{\rho_n^\sp} \sum_{\substack{\varphi\subseteq\PP_n\colon \#\varphi=k_0, \\ \dist(\varphi,\partial[0,1]^d) \le r_n}} \xi(r_n^{-1}\varphi) s'_n(\varphi,\PP_n) \\
&= \widetilde{H}_{n,k_0}^\sp - H_{n,k_0}^{\text{err},1}(\PP_n),
\end{split}
\end{align}
where we used that from \eqref{LOC} it follows that for all $\varphi\subseteq\PP_n$ with $\#\varphi=k_0$
\begin{equation}\label{equation_distances_sparse}
\xi(r_n^{-1}\varphi) t_n(\varphi,\PP_n) = \xi(r_n^{-1}\varphi) s'_n(\varphi,\PP_n).
\end{equation}
Further, we introduce the event
\begin{align*}
G_n := \big\{\#\{X\in\PP_n\setminus[r_n,1-r_n]^d\colon &s'_n(\{X\}\cup\varphi,\PP_n)=1 \text{ for some }\varphi\subseteq\PP_n \\
&\text{with } \#(\{X\}\cup\varphi)=k_0\} < \delta\rho_n^\sp\big\},
\end{align*}
which implies that the number of connected components of size $k_0$ with respect to the Euclidean distance that are located close to the boundary of $[0,1]^d$ is negligible. To deal with the probability of $G_n$, note that it is possible to replace the event in Lemma \ref{lemma_prob_resampling_between_boxes_sparse} with the complement of $G_n$ and we still get that
\begin{equation}\label{inequality_event_points_close_to_boundary_of_torus}
\limsup_{n\to\infty} \frac1{\rho_n^\sp} \log \P(G_n^c) = -\infty.
\end{equation}
To show this, the proof of Lemma \ref{lemma_prob_resampling_between_boxes_sparse} can be repeated with only one modification that arises from switching from the toroidal to the Euclidean metric. In \eqref{inequality_probability_bad_subbox_sparse} one has to consider that it is possible that only a fraction of the ball intersects the box.

Under $G_n$ the number of components summed over in $H_{n,k_0}^{\text{err},1}$ is bounded by $\delta\rho_n^\sp$ and with \eqref{inequality_fixed_size_error_sparse}, we can compute
\begin{align*}
\P(H_n^\sp \le a) &\le \P(\widetilde{H}_{n,k_0}^\sp - H_{n,k_0}^{\text{err},1}(\PP_n) \le a) \le \P(\widetilde{H}_{n,k_0}^\sp \le a + H_{n,k_0}^{\text{err},1}(\PP_n), G_n) + \P(G_n^c) \\
&\le \P\Big(\widetilde{H}_{n,k_0}^\sp \le a + \delta \sup_{\varphi\subseteq[0,1]^d, \#\varphi = k_0} \xi(r_n^{-1}\varphi), G_n\Big) + \P(G_n^c) \\
&\le \P(\widetilde{H}_{n,k_0}^\sp \le a + \delta b) + \P(G_n^c),
\end{align*}
where we used \eqref{INV}, \eqref{LOC} and \eqref{BND} to get for sufficiently large $n$
$$\sup_{\varphi\subseteq[0,1]^d, \#\varphi = k_0} \xi(r_n^{-1}\varphi) = \sup_{\varphi\subseteq[0,1]^d\setminus[r_n,1-r_n]^d, \#\varphi = k_0} \xi(r_n^{-1}\varphi) \le \sup_{\varphi\subseteq[0,1]^d, \#\varphi = k_0} \xi_n(\varphi) \le b.$$
By \eqref{inequality_event_points_close_to_boundary_of_torus}, the
probability of the complement of $G_n$ does not significantly contribute to the large deviations. From this point, \eqref{INV}, \eqref{LOC}, \eqref{BND} and \eqref{POS} let us apply \cite[Theorem 3.3]{hirschowada} to $\widetilde{H}_{n,k_0}^\sp$, which yields
$$\limsup_{n\to\infty} \frac1{\rho_n^\sp} \log \P(H_n^\sp \le a) \le \limsup_{n\to\infty} \frac1{\rho_n^\sp} \log \P(\widetilde{H}_{n,k_0}^\sp \le a + \delta b) \le -\inf_{\rho\colon T^\sp(\rho) \le a+\delta b} h^\sp(\rho\mid \tau_{k_0}^\sp).$$
and therefore, the asserted upper bound, after letting $\delta\rightarrow0$. Note that the rate function in \cite[Theorem 3.3]{hirschowada} is given as a Legendre transform. Arguing as in \cite[Corollary 3.2]{hirschowada}, this can be equivalently written in the relative entropy form.

For the lower bound, as a first step, with the same reasoning, we get for any $\delta>0$ that
\begin{equation}\label{inequality_LDP_sparse}
\liminf_{n\to\infty} \frac{1}{\rho_n^\sp} \log \P(\widetilde{H}_{n,k_0}^\sp < a-\delta) \geq -\inf_{\rho\colon T^\sp(\rho) < a-\delta} h^\sp(\rho \mid \tau_{k_0}^\sp).
\end{equation}
The next part of this proof is dedicated to show that in terms of large deviations, also for the lower bound, $H_n^\sp$ can be replaced with $\widetilde{H}_{n,k_0}^\sp$. For this, let
$$
E_n^\text{good} := \Big\{\max_{X\in \{Y\in Q\cap\PP'_n\colon Q\in\mathcal{J}_n\}}  \PP'_n\big(B_{2^d k_0 r_n}(X)\cap Q\big) \le k_0-1\Big\}
$$
and for $\e>0$ serving as parameter for the Bernoulli random variables,
$$
E_n := E_n^\text{good} \cap \bigcap_{Q\in \mathcal{J}_n} \{X_{Q, \eps} = 1\}\cap \bigcap_{Q\not\in \mathcal{J}_n} \{X_{Q, \eps} = 0\}.
$$
We start the computations with
$$\P(H_n^\sp < a) = \P(H_n^\sp(\PP''_n) < a) \ge \P(E_n, H_n^\sp(\PP''_n) < a).$$
Next, we can divide the functional into contributions that come from components intersecting the volume close to the boundary of a cube, denoted by
$$H_n^{\text{err},2}(\PP''_n) := \frac{1}{\rho_n^\sp} \sum_{\varphi\subseteq\PP''_n\colon \varphi \cap (\cup_{Q\in\mathcal{Q}_n} \partial_n Q) \neq \emptyset} \xi_n(\varphi) s_n(\varphi,\PP''_n),$$
and those that do not. Under the event $E_n$, we then have that
\begin{align*}
H_n^\sp(\PP''_n) &= \frac{1}{\rho_n^\sp} \sum_{\varphi\subseteq\PP''_n\colon \varphi \cap (\cup_{Q\in\mathcal{Q}_n} \partial_n Q) = \emptyset} \xi(r_n^{-1}\varphi) s_n(\varphi,\PP''_n) + H_n^{\text{err},2}(\PP''_n) \\
&\le \widetilde{H}_{n,k_0}^\sp(\PP_n) + H_n^{\text{err},2}(\PP''_n).
\end{align*}
We were able to bound the first term by $\widetilde{H}_{n,k_0}^\sp$ applied to $\PP_n$ instead of $\PP''_n$ because under the sprinkling event,  if we disregard the space close to the boundaries of the cubes, the coupled process $\PP''_n$ replaces $\PP_n$ in each cube that contained at least a part of a connected component of size $k_0+1$, without creating any new connected components of size $k_0$ or bigger. Further, we made use of \eqref{equation_distances_sparse} as in the proof of the upper bound. This lets us proceed with
$$\P(E_n, H_n^\sp(\PP''_n) < a) \ge \P(E_n, \widetilde{H}_{n,k_0}^\sp(\PP_n) + H_n^{\text{err},2}(\PP''_n) < a).$$
Further, to ease notation, let
$$F_n := \{\CCC_{n,k_0}(\cup_{Q\in \mathcal{Q}_n} \partial_n Q) < \delta\rho_n^\sp\}$$
denote the complement of the event from Lemma \ref{lemma_prob_resampling_between_boxes_sparse} for some $\delta>0$, which gives us
$$\P(E_n, H_{n,k_0}^\sp(\PP_n) + H_n^{\text{err},2}(\PP''_n) < a) \ge \P(E_n, F_n, \widetilde{H}_{n,k_0}^\sp(\PP_n) + H_n^{\text{err},2}(\PP''_n) < a).$$
Now, conditioned on $E_n$ and for sufficiently large $n$, the random geometric graph on $\PP''_n$ with connectivity radius $r_n$ cannot have a connected component of more than $2^d k_0$ nodes, since in that case if $n$ is large, a box $Q\in\mathcal{J}_n$ would exist that contains $k_0+1$ vertices of $\PP''_n$ with diameter less than or equal to $2^d k_0 r_n$. This contradicts $E_n^\text{good}$. Note that $2^d$ occurs here because it is the maximal number of boxes that can share a face. Thus, due to the nonnegativity of $\xi$, under $E_n\cap F_n$, it holds that
$$H_n^{\text{err},2}(\PP''_n)< \delta \sup_{\varphi\subseteq[0,1]^d, \#\varphi\le 2^d k_0} \xi_n(\varphi) \le \delta b,$$
where we recall that $b$ depending only on $d$ and $k_0$ arises from \eqref{BND}. This leads to
\begin{align*}
\P(E_n, F_n, \widetilde{H}_{n,k_0}^\sp(\PP_n) + H_n^{\text{err},2}(\PP''_n) < a) &\ge \P(E_n, F_n, \widetilde{H}_{n,k_0}^\sp(\PP_n) + \delta b < a) \\
&\ge \P(E_n, \widetilde{H}_{n,k_0}^\sp(\PP_n) + \delta b < a) - \P(F_n^c).
\end{align*}
Summarizing these steps and applying the tower property of the conditional expectation, we arrive at
$$
\P(H_n^\sp < a) \ge \E[\P(E_n | \PP_n) \one\{\widetilde{H}_{n,k_0}^\sp(\PP_n) < a-\delta b\}] - \P(F_n^c).$$
Now, using Lemma \ref{lemma_prob_new_edges_sparse} and independence of the events intersected in $E_n$ under $\PP_n$, we get that
$$
\P(E_n | \PP_n) = \P(E_n^\text{good} | \PP_n) \eps^{J_n} (1-\eps)^{\rho_n^\sp-J_n} \ge (\alpha_M \eps)^{J_n} (1-\eps)^{\rho_n^\sp}
$$
for an arbitrary $M > 2 \kappa_d^{k_0-1} 2^{k_0(d^2+1)} k_0^{k_0}$. This lets us proceed with
\begin{align*}
\P(H_n^\sp < a) &\ge \E\big[(\alpha_M \eps)^{J_n} (1-\eps)^{\rho_n^\sp} \one\{\widetilde{H}_{n,k_0}^\sp(\PP_n) < a-\delta b\}\big] - \P(F_n^c) \\
&\ge \E\big[(\alpha_M \eps)^{\delta\rho_n^\sp} (1-\eps)^{\rho_n^\sp} \one\{\widetilde{H}_{n,k_0}^\sp(\PP_n) < a-\delta b\} \one\{J_n < \delta\rho_n^\sp\}\big] - \P(F_n^c) \\
&\ge (\alpha_M \eps)^{\delta\rho_n^\sp} (1-\eps)^{\rho_n^\sp} \big(\P(\widetilde{H}_{n,k_0}^\sp(\PP_n) < a-\delta b) -\P(J_n \ge \delta\rho_n^\sp)\big) - \P(F_n^c).
\end{align*}
From this inequality and Lemma \ref{lemma_prob_resampling_between_boxes_sparse}, it follows that $\P(F_n^c)$ does not contribute significantly to the lower bound for the lower large deviations. Therefore, we arrive at
\begin{align*}
&\liminf_{n\to\infty} \frac{1}{\rho_n^\sp} \log \P(H_n^\sp < a) \\
\ge\ &\liminf_{n\to\infty} \frac{1}{\rho_n^\sp} \log\big((\alpha_M \eps)^{\delta\rho_n^\sp} (1-\eps)^{\rho_n^\sp} (\P(\widetilde{H}_{n,k_0}^\sp(\PP_n) < a-\delta b) -\P(J_n \ge \delta\rho_n^\sp))\big) \\
\ge\ &\delta \log (\alpha_M \eps) + \log(1-\eps) + \liminf_{n\to\infty} \frac{1}{\rho_n^\sp} \log\big(\P(\widetilde{H}_{n,k_0}^\sp(\PP_n) < a-\delta b) -\P(J_n \ge \delta\rho_n^\sp)\big).
\end{align*}
Now, Lemma \ref{lemma_probability_bad_boxes_sparse} implies that $\P(J_n \ge \delta\rho_n^\sp)$ does not affect the lower bound of the lower tails in this situation, and thus, plugging in \eqref{inequality_LDP_sparse}, we get
\begin{align*}
\liminf_{n\to\infty} \frac1{\rho_n^\sp} \log \P(H_n^\sp < a) &\ge \delta \log(\alpha_M \eps) + \log(1-\eps) + \liminf_{n\to\infty} \frac{1}{\rho_n^\sp} \log\big(\P(\widetilde{H}_{n,k_0}^\sp(\PP_n) < a-\delta b) \\
&\ge  \delta \log(\alpha_M \eps) + \log(1-\eps) - \inf_{\rho\colon T^\sp(\rho) < a-\delta b} h^\sp(\rho \mid \tau_{k_0}^\sp).
\end{align*}
Letting $\delta\rightarrow 0$ and then $\eps\rightarrow 0$ gives the lower bound
$$\liminf_{n\to\infty} \frac{1}{\rho_n^\sp} \log \P(H_n^\sp < a) \ge -\inf_{\rho\colon T^\sp(\rho) < a} h^\sp(\rho \mid \tau_{k_0}^\sp).$$
\enp

What follows are the proofs of the previously introduced lemmas. But, since we come across the task of bounding a similar quantity in the proofs of Lemmas \ref{lemma_probability_bad_boxes_sparse}, \ref{lemma_prob_new_edges_sparse} and \ref{lemma_prob_resampling_between_boxes_sparse}, we insert a short lemma that helps with this first.

\bel[Bound for the probability of many Poisson points in a ball]\label{lemma_bound_Poisson_points_in_ball}
For $Q\subseteq[0,1]^d$ and $m,l,r\in\N$, it holds that
$$\E[\#\{X\in Q\cap\PP_m\colon \PP_m(B_r(X)) \ge l\}] \le m^l \kappa_d^{l-1} r^{(l-1)d}  |Q|.$$
\enl

\bep[Proof of Lemma \ref{lemma_probability_bad_boxes_sparse}]
We are going to categorize boxes to create independence and use a binomial concentration inequality from \cite[Lemma 1.1]{poisson_conc}. We  use the set $\mathcal{L} := \{1,2\}^d$ to label each box in $\mathcal{Q}_n$ in a certain way to achieve that between two boxes of the same label, there will always be a box with a different label. To guarantee that this is possible on the torus, we assume that the number of boxes along each axis is divisible by $2$. For $l\in\mathcal{L}$, we denote the boxes of label $l$ by $\mathcal{Q}_n^{(l)}$. Then,
$$
	\P(J_n \ge \delta\rho_n^\sp) \le \sum_{l\in\mathcal{L}} \P(\#(\mathcal{Q}_n^{(l)}\cap \mathcal{J}_n) \ge \delta\rho_n^\sp/2^d).
$$
For $n$ large enough, the labeling guarantees that the events $\{Q\in\mathcal{J}_n\}$ are independent for different $Q\in\mathcal{Q}_n^{(l)}$. Thus, we are in a binomial setting and to use the mentioned binomial concentration inequality, we first bound the probability of one box being bad by using Lemma \ref{lemma_bound_Poisson_points_in_ball} to get that for an arbitrary $Q\in\mathcal{Q}_n$
\begin{align}\label{inequality_probability_bad_box_sparse}
\begin{split}
	\P(Q\in\mathcal{J}_n) &= \P(\PP_n(B_{2^d k_0 r_n}(X)) \ge k_0+1 \text{ for some } X\in Q\cap\PP_n) \\
	&\le \E\bigg[\sum_{X\in Q\cap\PP_n} \one\{\PP_n(B_{2^d k_0 r_n}(X)) \ge k_0+1\}\bigg] \\
	&\le n^{k_0+1} \kappa_d^{k_0}(2^d k_0 r_n)^{k_0d}  |Q| = \kappa_d^{k_0} (2^d k_0)^{k_0d} n r_n^d.
\end{split}
\end{align}
Next, using \cite[Lemma 1.1]{poisson_conc} for $n$ large, we get for every $l\in\mathcal{L}$ and every $\delta>0$, if $n$ is large enough, that
\begin{align*}
	\P(\#(\mathcal{Q}_n^{(l)}\cap \mathcal{J}_n) \ge \delta\rho_n^\sp/2^d) &\le \exp\bigg(-\frac{\delta\rho_n^\sp/2^d}{2} \log\Big(\frac{\delta\rho_n^\sp/2^d}{\rho_n^\sp \kappa_d^{k_0} (2^d k_0)^{k_0d} n r_n^d}\Big)\bigg) \\
	&= \exp\bigg(-\frac{\delta\rho_n^\sp}{2^{d+1}} \log\Big(\frac{\delta}{\kappa_d^{k_0} 2^d (2^d k_0)^{k_0d} n r_n^d}\Big)\bigg).
\end{align*}
The assumption $n r_n^d\rightarrow 0$ yields the assertion.
\enp

\bep[Proof of Lemma \ref{lemma_prob_new_edges_sparse}]
First, we let $M> 2 \kappa_d^{k_0-1} 2^{k_0(d^2+1)} k_0^{k_0}$ as well as $Q\in\mathcal{Q}_n$ be arbitrary and start by examining the probability that $\PP_{2n}$ has some amount of close points within $Q$ by invoking Markov's inequality and Lemma \ref{lemma_bound_Poisson_points_in_ball} to get
\begin{align*}
&\P\big(\#\{X\in Q\cap\PP_{2n} \colon \PP_{2n}(B_{2^d k_0 r_n}(X)\cap Q) \ge k_0\} \ge M\big) \\
&\le \frac1{M} \E\Big[\sum_{X\in Q\cap\PP_{2n}} \one\{\PP_{2n}(B_{2^d k_0 r_n}(X)\cap Q) \ge k_0\}\Big] \\
&\le \frac1{M} (2n)^{k_0} \kappa_d^{k_0-1} (2^d k_0 r_n)^{(k_0-1)d}  |Q| = \frac{\kappa_d^{k_0-1} 2^{k_0(d^2+1)} k_0^{k_0 d}}{M} \le \frac12.
\end{align*}
Note that a thinning of $\PP_{2n}$, where we keep each point independently with probability $1/2$ has the same distribution as $\PP_n$. Denote the thinned process by $\PP_{2n}^{\text{thin}}$. We proceed by deleting unwanted points in the thinning and get
\begin{align*}
&\P\Big(\max_{X\in Q\cap\PP'_n} \PP'_n(B_{2^d k_0 r_n}(X)\cap Q) \le k_0-1\Big) \\
&= \P\Big(\max_{X\in Q\cap\PP_{2n}^{\text{thin}}} \PP_{2n}^{\text{thin}}(B_{2^d k_0 r_n}(X)\cap Q) \le k_0-1\Big) \\
&= \E\big[(1/2)^{\#\{X\in Q\cap\PP_{2n} \colon \PP_{2n}(B_{2^d k_0 r_n}(X)\cap Q) \ge k_0\}}\big] \\
&\ge \E[2^{-M} \one\{\{X\in Q\cap\PP'_n \colon \PP'_n(B_{2^d k_0 r_n}(X)\cap Q) \ge k_0\} < M\}] \ge 2^{-M-1}.
\end{align*}
Now, we can use independence of the above events when considering different boxes to get
\begin{align*}
&\P\bigg(\bigcap_{Q\in\mathcal{J}_n} \Big\{\max_{X\in Q\cap\PP'_n} \PP'_n(B_{2^d k_0 r_n}(X)\cap Q) \le k_0-1\Big\} \biggm\vert \PP_n\bigg) \\
&= \prod_{Q\in\mathcal{J}_n} \P\Big(\max_{X\in Q\cap\PP'_n} \PP'_n(B_{2^d k_0 r_n}(X)\cap Q) \le k_0-1\Big) \ge (2^{-M-1})^{J_n}.
\end{align*}
\enp

\bep[Proof of Lemma \ref{lemma_prob_resampling_between_boxes_sparse}]
For a box $Q\in\mathcal{Q}_n$, we divide $\partial_n Q$ into a grid consisting of boxes with side length $r_n$ and call this collection of boxes $\mathcal{W}_n(\partial_n Q)$. We denote the total collection of these boxes by $\overline{\mathcal{W}}_n := \cup_{Q\in\mathcal{Q}_n} \mathcal{W}_n(\partial_n Q)$. Next, we can proceed with the same strategy that was already successfully applied in the proof of Lemma \ref{lemma_probability_bad_boxes_sparse}, but use more labels this time to achieve that between two boxes $W_1,W_2\in\overline{\mathcal{W}}_n$ of the same label, there are always $2^{d+1} k_0$ boxes labeled differently. We choose the label set $\mathcal{L} := \{1,2,\dots,2^{d+1} k_0 + 1\}^d$ and we reuse the notation $\overline{\mathcal{W}}_n^{(l)}$ for the boxes of label $l\in\mathcal{L}$. Again, we assume that the number of boxes along each axis is divisible by $2^{d+1}k_0+1$. This construction lets us search for connected components of at most $2^d k_0$ nodes in boxes with the same label independently. Our aim is to apply the already encountered binomial concentration bound \cite[Lemma 1.1]{poisson_conc} to the number of subcubes of a fixed label that contain a large connected component. This requires two things, a bound for the number of subcubes in $\overline{\mathcal{W}}_n^{(l)}$ and a bound for the probability of a subcube $W\in\overline{\mathcal{W}}_n$ containing at least one node in $\CCC_{n,k_0}$.

For the latter, i.e., the probability that $W\in\overline{\mathcal{W}}_n$ contains vertices that are part of a connected component of size between $k_0$ and $2^d k_0$, we compute, using Markov's inequality and Lemma \ref{lemma_bound_Poisson_points_in_ball}, that
\begin{align}\label{inequality_probability_bad_subbox_sparse}
\begin{split}
\P\big(W\cap\CCC_{n,k_0} \neq \emptyset\big) &\le \E[\#(W\cap\PP''_n\cap\CCC_{n,k_0})] \le \E[\#\{X\in W\cap\PP''_n\colon \PP''_n(B_{k_0 r_n}(X)) \ge k_0\}] \\
&\le n^{k_0} \kappa_d^{k_0-1} (k_0 r_n)^{d(k_0-1)} |W| = \kappa_d^{k_0-1} k_0^{(k_0-1)d} (n r_n^d)^{k_0}.
\end{split}
\end{align}
To find a bound for the number of subcubes, note that the volume of $\partial_n Q$ for $Q\in\mathcal{Q}_n$ is of order
$$(\rho_n^\sp)^{-(d-1)/d} r_n = n^{-k_0(d-1)/d} r_n^{-(k_0-1)(d-1)+1}.$$
Consequently, the number of boxes in $\mathcal{W}_n(\partial_n Q)$ can be bounded by dividing the above by the volume of a subcube $r_n^d$, which yields
\begin{equation}\label{inequality_number_bad_subboxes_sparse}
\#\mathcal{W}_n(\partial_n Q) \le c_1 n^{-k_0(d-1)/d} r_n^{-k_0(d-1)} = c_1 (n r_n^d)^{-k_0(d-1)/d}
\end{equation}
for some $c_1 := c_1(d, k_0) > 0$.  Thus, there are at most $\rho_n^\sp c_1 (n r_n^d)^{-k_0(d-1)/d}$ subcubes in $\overline{\mathcal{W}}_n$.

Before we invoke \cite[Lemma 1.1]{poisson_conc}, we can union over all labels and combine this with the union bound to get
\begin{align*}
	\P\big(\CCC_{n,k_0}(\cup_{Q\in \mathcal{J}_n} \partial_n Q) \ge \delta\rho_n^\sp \big)
	&\le \P\bigg(\bigcup_{l\in\mathcal{L}} \big\{\CCC_{n,k_0}(\cup_{W\in \overline{\mathcal{W}}_n^{(l)}} W) \ge \delta\rho_n^\sp/(\#\mc L)\big\}\bigg) \\
	&\le \sum_{l\in\mathcal{L}} \P\big(\CCC_{n,k_0}(\cup_{W\in \overline{\mathcal{W}}_n^{(l)}} W) \ge \delta\rho_n^\sp/(\#\mc L)\big).
\end{align*}
At this point,  let $W\in\overline{\mathcal{W}}_n$ be arbitrary. An important observation is that $\CCC_{n,k_0}(W)$ is bounded by a constant that does not depend on $n$. More precisely,  a connected component occupies a ball of radius at least $r_n$ that cannot intersect any other connected component. Consequently, when choosing $r_n/2$ as radius instead, that ball cannot intersect any ball of radius $r_n/2$ that is centered at a node that belongs to another connected component. When considering connected components with a vertex in $W$, at least $1/2^d$ of the volume of a ball with radius $r_n/2$ centered at that vertex has to be contained in $W$. The factor $1/2^d$ adjusts for the possibility that the center of the ball is in a corner of $W$. Therefore, we can bound the available space by $|W|$ and the maximal component size by $2^d k_0$ and arrive at
$$\CCC_{n,k_0}(W) \le \frac{2^d k_0 |W|}{\kappa_d (r_n/2)^d /2^d} = 8^d \kappa_d^{-1} k_0,$$
which implies that for a fixed $l\in\mathcal{L}$
$$
\P\big(\CCC_{n,k_0}(\cup_{W\in \overline{\mathcal{W}}_n^{(l)}} W) \ge \delta\rho_n^\sp/(\#\mc L)\big) \le \P\big(\#\{W\in \overline{\mathcal{W}}_n^{(l)} \colon W\cap\CCC_{n,k_0} \neq \emptyset\} \ge \delta\rho_n^\sp/c_2 \big),
$$
where $c_2 := c_2(d,k_0) := 8^d\kappa_d^{-1} k_0 \#\mc L$.
Next, the independence guaranteed by the labeling and the bounds derived in \eqref{inequality_probability_bad_subbox_sparse} and \eqref{inequality_number_bad_subboxes_sparse} let us apply the binomial bound \cite[Lemma 1.1]{poisson_conc} for sufficiently large $n$ to arrive at
\begin{align*}
	&\P\big(\#\{W\in \overline{\mathcal{W}}_n^{(l)} \colon W\cap\CCC_{n,k_0} \neq \emptyset\} \ge \delta\rho_n^\sp/c_2\big) \\
	\le\ &\exp\bigg(-\frac{\delta\rho_n^\sp}{2 c_2} \log\Big(\frac{\delta\rho_n^\sp/c_1}{\rho_n^\sp d c_1 (n r_n^d)^{-k_0(d-1)/d} \kappa_d^{k_0-1} k_0^{(k_0-1)d} (n r_n^d)^{k_0}}\Big)\bigg) \\
	=\ &\exp\bigg(-\frac{\delta\rho_n^\sp}{2 c_2} \log\Big(\frac{\delta}{c_2 c_1 \kappa_d^{k_0-1} k_0^{(k_0-1)d} (n r_n^d)^{k_0/d}}\Big)\bigg),
\end{align*}
yielding the assertion, since $n r_n^d \rightarrow 0$.
\enp

\bep[Proof of Lemma \ref{lemma_bound_Poisson_points_in_ball}]
If $l=1$, we get
\begin{align*}
&\E[\#\{X\in Q\cap\PP_m\colon\PP_m(B_r(X)) \ge 1\}] = \E[\PP_m(Q)] = m |Q|.
\end{align*}
For $l>1$, an application of Mecke's equation and Markov's inequality yields
\begin{align*}
&\E[\#\{X\in Q\cap\PP_m\colon\PP_m(B_r(X)) \ge l\}] = m \int_{Q} \E[\one\{\PP_m(B_r(x)) \ge l-1\}] \d x \\
	&= m \int_{Q} \P(\{Y_1,\dots,Y_{l-1}\}\subseteq B_r(x) \text{ for some } \{Y_1,\dots,Y_{l-1}\}\subseteq\PP_m) \d x \\
	&\le m \int_Q \E\bigg[\sum_{\{Y_1,\dots,Y_{l-1}\}\subseteq \PP_m} \one\{Y_1,\dots,Y_{l-1}\in B_r(x)\}\bigg] \d x \\
	&\le m^l \int_{Q} \int_{[0,1]^{(l-1)d}} \one\{y_1,\dots,y_{l-1}\in B_r(x)\} \d (y_1,\dots,y_{l-1}) \d x = m^l \kappa_d^{l-1} r^{(l-1)d}  |Q|.
\end{align*}
\enp
\section{Proof of Theorem \ref{theorem_main_dense} (dense)}\label{section_proof_dense}

The general outline of the proof of the dense regime follows the ideas for the sparse case. Here, we aim to apply a contraction principle using the large deviation asymptotics with respect to the weak topology from \cite{hirschowadakang}. In order to apply the contraction principle directly, $T_k^\de$ must be continuous with respect to the weak topology, meaning that the integrand needs to be bounded. However, this condition is not immediately satisfied. To overcome this, using the technique of sprinkling, we aim to artificially introduce a bound for the score function that will translate to the integrand of $T_k^\de$.

As in the sparse regime, we divide $[0,1]^d$ into a grid of cubes with side length $(\rho_{n,k}^\de)^{-1/d}$, assuming that $\rho_{n,k}^\de$ is a natural number and denote this collection by $\mathcal{Q}_n$. We are going to use the same objects that were introduced in the sparse regime. As a reminder, $\PP'_n$ is a Poisson point process on $[0,1]^d$ with intensity $n$ independent of $\PP_n$, and for Bernoulli random variables with parameter $\eps\in(0,1)$, independent of each other and all introduced Poisson random measures, for every $Q\in\mc Q_n$, we defined
$$\PP_n^{Q} := \begin{cases}
	Q\cap\PP'_n &\text{if } X_{Q, \eps} = 1 \\
	Q\cap\PP_n &\text{if } X_{Q, \eps} = 0
	\end{cases}.$$
Finally, we denoted the union $\cup_{Q\in\mathcal{Q}_n} \PP_n^Q$ by $\PP''_n$. In the dense regime, we aim to use the Bernoulli random variables to control $\PP''_n$ in such a way that we resample $\PP_n$ using $\PP'_n$ in each box that makes it too likely that there is an $X$ with a large edge while keeping $\PP_n$ in all other boxes. Mathematically expressed, for a random configuration $\eta\in\NNN$, we want to avoid boxes that foster the existence of an $X\in\eta$ with
\begin{equation}\label{equation_indicator}
\xi_n(X, \eta): = (n\kappa_dR_k(X, \eta)^d - a_n - s_0)_+ > M
\end{equation}
for $M>0$. If a box $Q\in \mc Q_n$ has no such point within $Q\cap\eta$, we will refer to it as $(\eta,M)$\emph{-bounded}.
To achieve this goal, we need to ensure that the resampling is done in such a way that adjacent boxes remain compatible in the sense that even after the resampling, the conditional probability that a box fulfills the boundedness property remains high. To that end, we fix an arbitrary ordering of the boxes in $\mc Q_n$ such that $Q_n^{(i)}$ denotes the $i$th box in $\mc Q_n$ and then impose conditions recursively. More precisely, for $\eta\in\{\PP_n,\PP'_n\}$, we denote $\bar{\mc S}_n^{(1)}(\eta) := \mc (Q_n^{(1)}\cap\eta)\cup(\PP_n\setminus Q_n^{(1)})$, where outside of the box $Q_n^{(1)}$ we could have used an arbitrary Poisson point process with intensity $n$ in the definition of $\bar{\mc S}_n^{(1)}$. Next, for an arbitrary $j\in\{1,\dots,\rho_{n,k}^\de\}$, let
\begin{enumerate}
\item $\mc N(j)$ denotes the ordering indices of the boxes adjacent to box $Q_n^{(j)}$;
\item $\mc N_+(j) := \mc N(j)\cup\{j\}$ be the above unioned with $\{j\}$;
\item $\eta^{(j)} := \cup_{s\le j} (Q_n^{(s)}\cap\eta)$ be $\eta$ restricted to the first $j$ boxes.
\end{enumerate}
Then, by setting
$$a_{i, 1}(\eta) := \P\big(Q_n^{(i)}\text{ is $(\bar{\mc S}_n^{(1)}(\eta),M)$-bounded}\mid \PP_n^{(1)}, (\PP'_n)^{(1)}\big), \quad i\in\NN_+(1),$$
we label the box $Q_n^{(1)}$ as \emph{$(\eta,M)$-good} if 
\begin{align*}
	\min_{i \in \NN_+(1)}a_{i, 1}(\eta) \ge 1 - e^{-M/2}.
\end{align*}
Then, we proceed step by step and for $2\le j\le \rho_{n,k}^\de$ set
$$\mc S_n^{(j-1)} := \begin{cases}
	Q_n^{(j-1)}\cap\PP'_n \quad\text{if } Q_n^{(j-1)} \text{ is }(\PP_n,M)\text{-bad} \\
	Q_n^{(j-1)}\cap\PP_n \quad\text{if } Q_n^{(j-1)} \text{ is }(\PP_n,M)\text{-good}
	\end{cases},$$
to be able to define
$$
\bar{\mc S}_n^{(j)}(\eta) := (Q_n^{(j)}\cap\eta) \cup (\cup_{s\le j-1} \mc S_n^{(s)}) \cup (\PP_n\setminus \cup_{s\le j} Q_n^{(j)}).
$$
Additionally, for $i\in\NN_+(j)$, we define the conditional probabilities
$$a_{i, j}(\eta) := \P\big(Q_n^{(i)}\text{ is $(\bar{\mc S}_n^{(j)}(\eta),M)$-bounded}\mid \PP_n^{(j)}, (\PP'_n)^{(j)}\big)$$
and note that $a_{i, j}(\eta)$ only depends on the configurations of $\PP_n,\PP'_n$ in $Q_n^{(s)}$ for $s\in\NN_+(i)\cap\{1,\dots,j\}$. We then say that the box $Q_n^{(j)}$ is \emph{$(\mc \eta, M)$-good} if 
\begin{align}
	\label{eq:mgood}
	\min_{i \in \NN_+(j)} a_{i, j}(\eta) \ge 1 - b_{i,j}^{(M)}
\end{align}
holds, where
\begin{equation}\label{equation_probability_Mgood}
b_{i,j}^{(M)} := e^{-M 2^{-1-\#\{s\in\NN_+(i) \colon s\le j\}}}.
\end{equation}
In words, we consider a configuration $\eta$ within the box $Q_n^{(j)}$ suitable if the probability of any adjacent box $Q_n^{(i)}$ being $(\bar{\mc S}_n^{(j)}(\eta),M)$-bounded is large conditioned on the configurations in the boxes that have already been considered in a step $s < j$ and the configuration $Q_n^{(j)}\cap\eta$ in the current box. 

Next, let 
$$\mathcal{J}_n^M := \mathcal{J}_n^M(\PP_n, \PP_n') := \{Q_n^{(i)} \in \mathcal{Q}_n\colon\text{ $Q_n^{(i)}$ is $(\PP_n,M)$-bad}\}$$ 
be the collection of $(\PP_n,M)$-bad boxes and we abbreviate its cardinality by
$$J_n^M := \#\mathcal{J}_n^M.$$
Since $k$ can be considered as fixed now, we can write
$$\rho_n^\de := \rho_{n,k}^\de$$
to ease notation.
We first make sure that those bad boxes do not occur too many times with a probability that is too high.

\bel[Bad boxes are exponentially negligible]\label{lemma_probability_bad_boxes}
Let $\delta, M>0.$ Then,
$$\P(J_n^M \ge \delta\rho_n^\de) \le \exp\Big(-\frac{\delta\rho_n^\de}{5^d 2} \log\Big(\frac{\delta e^{M/2+s_0}}{15^d 2^{3^d+1} k}\Big)\Big).$$
In particular,
$$\limsup_{M\to\infty}\limsup_{n\to\infty} \frac1{\rho_n^\de} \log\P(J_n^M \ge \delta\rho_n^\de) = -\infty.$$
\enl

Furthermore, we do not desire that a resampled box is still deemed bad. To achieve this, for $j\in\{1,\dots,\rho_n^\de\}$ and $\eta$ equal to either $\PP_n$ or $\PP'_n$, we let
\begin{equation}\label{equation_good_condition1}
E_j^\text{good}(\eta) := E_{j,n}^\text{good}(\eta) :=  \{Q_n^{(j)}\text{ is }(\eta,M)\text{-good}\}
\end{equation}
be the event that $Q_n^{(j)}$ is $(\eta,M)$-good. In addition, let $E_j^\text{bad}(\eta)$ be the event's complement and for $M_0>0$, let
\begin{equation}\label{equation_good_condition2}
E_j^{b} := E_{j,n}^{b} := \{Q_n^{(j)}\setminus\partial_n Q_n^{(j)}\text{ is }(\PP'_n, M_0)\text{-bounded}\}
\end{equation}
be the event that $\PP'_n$ not close to the boundary of a box $Q_n^{(j)}$ fulfills an additional boundedness condition. Here, for every $Q\in\mathcal{Q}_n$, we denoted by
$$\partial_n Q := \big\{x\in Q\colon \dist(\{x\},\partial Q) \le t_n \big\}$$
the set of all points in $Q$ within distance 
$$t_n := \Big(\frac{a_n + w_n}{n \kappa_d }\Big)^{1/d}$$
of the complement of $Q$, where $(w_n)_n$ is a sequence with $w_n \rightarrow \infty$ and $w_n \in o(a_n)$ that we henceforth fix.

The next lemma states that for each $n\in\N$, conditioned on $\PP_n$, the probability that a box is either good, or we can resample it in a beneficial way otherwise is positive. 

\bel[Lower bound for probability of a good box or resampling a good box]\label{lemma_prob_new_edges}
For any $M , M_0 > 0$ it holds that
\begin{align*}
\P\bigg(\bigcap_{j=1}^{\rho_n^\de} E_j^\text{good}(\PP_n) \cup \big(E_j^\text{bad}(\PP_n) \cap E_j^\text{good}(\PP'_n) \cap E_j^b \big) \biggm\vert \PP_n\bigg) \ge q_{M_0,M}^{\rho_n^\de},
\end{align*}
where $q_{M_0,M} := 1 - 3^d 2 k e^{|s_0|} (e^{-M_0} + e^{-M/2^{4^d}})$.
\enl

Then, for $\e>0$ serving as parameter for the Bernoulli random variables, we define
$$
E_n^* := \bigcap_{i=1}^{\rho_n^\de}\big(E_j^\text{good}(\PP_n) \cap \{X_{Q_n^{(j)},\e} = 0\}\big) \cup \big(E_j^\text{bad}(\PP_n) \cap E_j^\text{good}(\PP'_n) \cap E_j^b \cap \{X_{Q_n^{(j)},\e} = 1\}\big).
$$
Recalling the definition of the mixed Poisson point process $\PP''_n$, this means that, using the Bernoulli random variables, we resample all boxes that are bad with respect to $\PP_n$ and ask for $\PP'_n$ to satisfy the goodness as in the event in \eqref{equation_good_condition1} and the additional condition described in \eqref{equation_good_condition2} in the boxes, where the sprinkling triggered.

\bel[Lower bound for probability of resampling bad boxes]\label{lemma_lower_bound_good_event}
For $\e\in(0,1)$ and arbitrary $M,M_0>0$, the event $E_n^*$ satisfies that
\begin{equation}\label{subset_good_event}
E_n^* \subseteq \{Q_n^{(i)} \text{ is } (\PP''_n,M)\text{-bounded for every }i\le\rho_n^\de\}.
\end{equation}
Further, it holds that
\begin{equation}\label{inequality_good_event}
\P(E_n^* \mid \PP_n) \ge \e^{\delta\rho_n^\de} (1-\e)^{\rho_n^\de} \big(q_{M_0,M}^{\rho_n^\de} - \P(J_n^M \ge \delta\rho_n^\de \mid \PP_n)\big).
\end{equation}
\enl

Recalling the definition of $\xi_n$ in \eqref{equation_indicator}, we introduce the error terms
\begin{equation}\label{equation_dense_error1}
	H_{n,M}^{\text{err},\partial}(\PP''_n) := \frac1{\rho_n^\de}\sum_{ X \in \PP_n''  \cap (\cup_{Q\in\mathcal{Q}_n} \partial_n Q) } M \wedge \xi_n(X,\PP''_n)  
\end{equation}
and 
\begin{equation}\label{equation_dense_error2}
	H_{n,M,M_0}^{\text{err}, \mc J}(\PP_n, \PP'_n, \PP''_n) := \frac1{\rho_n^\de}\sum_{ X \in \PP_n''  \cap (\cup_{Q\in\mathcal{J}_n^M} Q ) } M_0 \wedge \xi_n(X,\PP''_n) 
\end{equation}
that will denote potential deviations introduced by the sprinkling. The following lemma is devoted to show that these errors are insignificant.

\bel[$H_{n,M}^{\text{err},\partial}(\PP''_n)$ and $H_{n,M,M_0}^{\text{err}, \mc J}(\PP_n, \PP'_n, \PP''_n)$ are negligible]\label{lemma_prob_resampling_between_boxes}
Let $\delta>0$. Then, for any $M>0$
$$\limsup_{n\to\infty} \frac1{\rho_n^\de} \log \P\big(H_{n,M}^{\text{err},\partial}(\PP''_n) \ge \delta\big) = -\infty$$
and for additionally any $M_0>0$
$$
\limsup_{M\to\infty}\limsup_{n\to\infty} \frac1{\rho_n^\de} \log \P\big(H_{n,M,M_0}^{\text{err}, \mc J}(\PP_n, \PP'_n, \PP''_n) \ge \delta\big) = -\infty.
$$
\enl

These lemmas allow us to prove the main theorem.

\bep[Proof of Theorem \ref{theorem_main_dense}]
Let $M>0$. We start by defining the functional
$$H_{n,M} := H_{n,M}(\PP_n) :=\frac{1}{\rho_n^\de} \sum_{X \in \PP_n} \xi_n(X, \PP_n) \one\{n\kappa_d R_k(X, \PP_n)^d - a_n -s_0 \le M\},$$
where we only add up scores of vertices, for which the distance to the $k$-closest node satisfies an additional bound, with the goal of applying \cite[Theorem 2.1]{hirschowada} to it. Along these lines, we define
$$L_{n,k} := \frac1{\rho_n^\de} \sum_{X\in\PP_n} \delta_{n\kappa_d R_k(X,\PP_n)^d - a_n}$$
as a random Radon measure on $\R$, which we henceforth restrict to a random Radon measure on $E_0$, denoted by $L_{n,k}^{E_0}$. Next, defined on the domain of Radon measures on $E_0$, the map given by
$$T_M(\rho) := \int_{E_0} (x-s_0) \wedge M \d\rho(x)$$
is continuous with respect to the weak topology and applied to $L_{n,k}^{E_0}$ yields $T_M(L_{n,k}^{E_0}) = H_{n,M}$.

Now, for the upper bound, note that
$$H_{n,M} \le H_n^\de.$$
From this point, \cite[Theorem 2.1]{hirschowada} and the contraction principle yield
$$\limsup_{n\to\infty} \frac{1}{\rho_n^\de} \log \P(H_n^\de \le a) \le \limsup_{n\to\infty} \frac{1}{\rho_n^\de} \log \P(H_{n,M} \le a) \le -\inf_{\rho\colon T_M(\rho) \le a} h^\de(\rho \mid \tau_k^\de)$$
and therefore,
$$\limsup_{n\to\infty} \frac{1}{\rho_n^\de} \log \P(H_n^\de \le a) \le - \limsup_{M\to\infty} \inf_{\rho\colon T_M(\rho) \le a} h^\de(\rho \mid \tau_k^\de).$$
Using monotone convergence of $T_M(\rho)$ towards $T_k^\de(\rho)$ for every Radon measure $\rho$ on $E_0$ as $M\rightarrow\infty$, gives the assertion.

For the lower bound, with the same reasoning we get for any $\delta>0$ that
\begin{equation}\label{inequality_LDP}
\liminf_{n\to\infty} \frac{1}{\rho_n^\de} \log \P(H_{n,M} < a-\delta) \geq -\inf_{\rho\colon T_M(\rho) < a-\delta} h^\de(\rho \mid \tau_k^\de).
\end{equation}
Next, as in the proof of the sparse regime, we need to show that $H_n^\de$ can be replaced with $H_{n,M}$ when it comes to the lower large deviations. 
We start the computations with
$$\P(H_n^\de < a) = \P(H_n^\de(\PP''_n) < a) \ge \P(E_n^*, H_n^\de(\PP''_n) < a).$$
Next, let $M_0>0$. Then, under the event $E_n^*$ we assert that
\begin{equation}\label{inequality_good_event_errors}
	H_n^\de(\PP''_n)  \le H_{n,M}(\PP_n) + H_{n,M}^{\text{err},\partial}(\PP''_n) + H_{n,M,M_0}^{\text{err}, \mc J}(\PP_n, \PP'_n, \PP''_n),
\end{equation}
where we recall the definitions of the error terms $H_{n,M}^{\text{err},\partial}(\PP''_n)$ and $H_{n,M,M_0}^{\text{err}, \mc J}(\PP_n, \PP'_n, \PP''_n)$ from \eqref{equation_dense_error1} and \eqref{equation_dense_error2}.
To show this claim, we partition $[0,1]^d$ into three subsets. Let
\begin{enumerate}
\item $S_1 := \cup_{Q\in\mc Q_n} \partial_n Q$, be the space close to the boundary of each box;
\item $S_2 := \cup_{Q\in\mathcal{J}_n^M} Q\setminus\partial_n Q$, be the union of all bad boxes without the space close to their boundaries;
\item $S_3 := \cup_{Q\in\mc Q_n\setminus\mathcal{J}_n^M} Q\setminus\partial_n Q$, be the union of all good boxes without the space close to their boundaries.
\end{enumerate}
Then,
$$
H_n^\de(\PP''_n) = \underbrace{ \frac1{\rho_{n, k}^\de}\sum_{ X \in \PP_n''  \cap S_1} \xi_n(X,\PP''_n)}_{=: (\star)} + \underbrace{ \frac1{\rho_{n, k}^\de}\sum_{ X \in \PP_n''  \cap S_2} \xi_n(X,\PP''_n)}_{=: (\star\star)} + \underbrace{ \frac1{\rho_{n, k}^\de}\sum_{ X \in \PP_n''  \cap S_3} \xi_n(X,\PP''_n)}_{=: (\star\star\star)}.
$$
Under $E_n^*$, for all $X\in\PP''_n\cap S_1$ it is satisfied that the box in which $X$ is located is $(\PP''_n,M)$-bounded by Lemma \ref{lemma_lower_bound_good_event}, which means that $\xi_n(X,\PP''_n) \le M$. Thus,
$$(\star) \le H_{n,M}^{\text{err},\partial}(\PP''_n).$$
Further, for all boxes $Q\in\mc Q_n\setminus\mathcal{J}_n^M$, i.e., that are already $(\PP_n,M)$-good, we stress that the distance of $\partial_n Q$ to the boundary of $Q$ was set to be at least $t_n$, and thus, we can assume that this distance is larger than $((M+a_n+s_0)/(n\kappa_d))^{1/d}$. Therefore, points in $\partial_n Q$ for a $(\PP_n,M)$-good box $Q$ are not affected by the potential replacement of $\PP_n$ with $\PP'_n$ in adjacent boxes, which means that due to the $(\PP_n,M)$-boundedness of $Q$, all nodes $X\in (Q\setminus\partial_n Q)\cap\PP_n$ satisfy that $\xi_n(X,\PP_n)\le M$. This yields that for large enough $n$
$$(\star\star\star) \le H_{n,M}(\PP_n).$$
Finally, under $E_n^*$, for all boxes $Q$ that were initially $(\PP_n,M)$-bad, the sprinkling assures that $Q\setminus\partial_n Q$ is $(\PP''_n,M_0)$-bounded, which results in
$$(\star\star) \le H_{n,M,M_0}^{\text{err}, \mc J}(\PP_n, \PP'_n, \PP''_n)$$
and confirms \eqref{inequality_good_event_errors}.

This lets us proceed with
$$\P(E_n^*, H_n^\de(\PP''_n) < a) \ge \P(E_n^*, H_{n,M}(\PP_n) + H_{n,M}^{\text{err},\partial}(\PP''_n) + H_{n,M,M_0}^{\text{err}, \mc J}(\PP_n, \PP'_n, \PP''_n) < a).$$
Further, to ease notation, let
$$F_n := \{H_{n,M}^{\text{err},\partial}(\PP''_n)< \delta\} \cap \{H_{n,M,M_0}^{\text{err}, \mc J}(\PP_n, \PP'_n, \PP''_n) < \delta\}$$
denote the complements of the events from Lemma \ref{lemma_prob_resampling_between_boxes} for some $\delta>0$, which gives us
\begin{align*}
&\P(E_n^*, H_{n,M}(\PP_n) + H_{n,M}^{\text{err},\partial}(\PP''_n) + H_{n,M,M_0}^{\text{err}, \mc J}(\PP_n, \PP'_n, \PP''_n) < a) \\
&\ge \P(E_n^*, F_n, H_{n,M}(\PP_n) + 2\delta < a) \\
&\ge \P(E_n^*, H_{n,M}(\PP_n) < a - 2\delta) - \P(F_n^c).
\end{align*}
Summarizing these steps and applying the tower property of the conditional expectation, we arrive at
$$
\P(H_n^\de < a) \ge \E[\P(E_n^* \mid \PP_n) \one\{H_{n,M}(\PP_n) < a-2\delta\}] - \P(F_n^c).$$
Now, due to Lemma \ref{lemma_lower_bound_good_event}, we get that
\begin{align}\label{inequality_ldp_with_errors}
\begin{split}
&\P(H_n^\de < a) \\
&\ge \e^{\delta\rho_n^\de} (1-\e)^{\rho_n^\de} \Big(q_{M_0,M}^{\rho_n^\de} \P(H_{n,M}(\PP_n) < a-2\delta) \\
&\qquad\qquad\qquad\qquad\quad - \E[\P(J_n^M \ge \delta\rho_n^\de \mid \PP_n) \one\{H_{n,M}(\PP_n) < a-2\delta\}]\Big) - \P(F_n^c) \\
&\ge \e^{\delta\rho_n^\de} (1-\e)^{\rho_n^\de} q_{M_0,M}^{\rho_n^\de} \P(H_{n,M}(\PP_n) < a-2\delta)  - \P(J_n^M \ge \delta\rho_n^\de) - \P(F_n^c).
\end{split}
\end{align}
From here, Lemmas \ref{lemma_probability_bad_boxes} and \ref{lemma_prob_resampling_between_boxes} assert that neither $\P(J_n^M \ge \delta\rho_n^\de)$ nor $\P(F_n^c)$ contribute significantly to the lower bound for the lower large deviations. Thus, we focus on the first term of the sum in the last line of \eqref{inequality_ldp_with_errors} and examine it under the assumption that $M$ and $M_0$ are large enough such that $q_{M_0,M}>0$ by computing
\begin{align*}
&\liminf_{M\to\infty} \liminf_{n\to\infty} \frac1{\rho_n^\de} \log\big(\e^{\delta\rho_n^\de} (1-\e)^{\rho_n^\de} q_{M_0,M}^{\rho_n^\de} \P(H_{n,M}(\PP_n) < a-2\delta)\big) \\
&\ge \delta \log \e + \log(1-\e) + \log q_{M_0,\infty} + \liminf_{M\to\infty} \liminf_{n\to\infty} \frac1{\rho_n^\de} \log \P(H_{n,M}(\PP_n) < a-2\delta),
\end{align*}
where $q_{M_0,\infty} := 1 - 3^d 2 k e^{|s_0|} e^{-M_0}$.
Now, after plugging in \eqref{inequality_LDP}, we arrive at
\begin{align*}
&\liminf_{M\to\infty} \liminf_{n\to\infty} \frac1{\rho_n^\de} \log\big(\e^{\delta\rho_n^\de} (1-\e)^{\rho_n^\de} q_{M_0,M}^{\rho_n^\de} \P(H_{n,M}(\PP_n) < a-2\delta)\big) \\
&\ge \delta \log \e + \log(1-\e) + \log q_{M_0,\infty} - \limsup_{M\to\infty} \inf_{\rho\colon T_M(\rho) < a-2\delta} h^\de(\rho \mid \tau_k^\de) \\
&\ge \delta \log \e + \log(1-\e) + \log q_{M_0,\infty} - \inf_{\rho\colon T_k^\de(\rho) < a-2\delta} h^\de(\rho \mid \tau_k^\de),
\end{align*}
where in the last line we used that $T_M(\rho)\le T_k^\de(\rho)$.
Letting $\delta\rightarrow 0$, $\eps\rightarrow 0$ and then $M_0\rightarrow \infty$ gives the lower bound
$$\liminf_{n\to\infty} \frac{1}{\rho_n^\de} \log \P(H_n^\de < a) \ge -\inf_{\rho\colon T_k^\de(\rho) < a} h^\de(\rho \mid \tau_k^\de).$$
\enp

What follows are the proofs of the previously introduced lemmas.
\bep[Proof of Lemma \ref{lemma_probability_bad_boxes}]
We claim that for some $c := c(d,k) >0$
\begin{equation}\label{inequality_prob_bad_box}
\P(Q^{(j)} \text{ is } (\PP_n,M)\text{-bad}) \le c e^{-M/2-s_0}
\end{equation}
if we choose $n$ sufficiently large. Once the claim in \eqref{inequality_prob_bad_box} is established, we conclude the proof as follows. For each $n\in\N$, we will categorize the boxes in $\mc Q_n$ to create independence and use the already encountered binomial concentration inequality from \cite[Lemma 1.1]{poisson_conc}. We can use $5^d$ labels, for instance, the set $\mathcal{L} := \{1,2,3,4,5\}^d$, to label each box in $\mathcal{Q}_n$ in a certain way to achieve that between two boxes of the same label, there will always be four boxes with different labels. Here, we assumed that the number of boxes along each axis is divisible by $5$. For $l\in\mathcal{L}$, we denote the boxes of label $l$ by $\mathcal{Q}_n^{(l)}$. Then,
$$
	\P(J_n^M \ge \delta\rho_n^\de) \le \sum_{l\in\mathcal{L}} \P(\#(\mathcal{Q}_n^{(l)}\cap \mathcal{J}_n^M) \ge \delta\rho_n^\de/5^d).
$$
For $n$ large enough, the labeling guarantees that the events $\{Q\in\mathcal{J}_n^M\}$ are independent for different $Q\in\mathcal{Q}_n^{(l)}$. Thus, we are in a binomial setting and can invoke \cite[Lemma 1.1]{poisson_conc} with success probability given by the bound in \eqref{inequality_prob_bad_box}, to get for every $l\in\mathcal{L}$ and delta $\delta>0$ that
\begin{align*}
	\P(\#(\mathcal{Q}_n^{(l)}\cap \mathcal{J}_n^M) \ge \delta\rho_n^\de/5^d) &\le \exp\bigg(-\frac{\delta\rho_n^\de/5^d}{2} \log\Big(\frac{\delta\rho_n^\de/5^d}{\rho_n^\de  ce^{-M/2-s_0}}\Big)\bigg) \\
	&= \exp\bigg(-\frac{\delta\rho_n^\de}{5^d 2} \log\Big(\frac{\delta e^{M/2+s_0}}{5^d c}\Big)\bigg)
\end{align*}
if $n$ is large enough. From this point, we see that
$$\frac1{\rho_n^\de} \log\P(J_n^M \ge \delta\rho_n^\de) \le -\frac{\delta}{5^d 2} \log\Big(\frac{\delta e^{M/2+s_0}}{5^d c}\Big),$$
and the right-hand side does not depend on $n$ anymore.  Furthermore, it satisfies that
$$-\frac{\delta}{5^d 2} \log\Big(\frac{\delta e^{M/2+s_0}}{5^d c}\Big) \overset{M\to\infty}{\longrightarrow} -\infty.$$
It remains to show \eqref{inequality_prob_bad_box}. For this, let $j\le \rho_n^\de$ and $i\in\NN_+(j)$ be arbitrary. Then, the tower property yields
\begin{align*}
&\P(Q_n^{(i)} \text{ is }(\bar{\mc S}_n^{(j)}(\PP_n),M)\text{-bounded}) = \E[a_{i,j}(\PP_n)] \\
&= \E[a_{i,j}(\PP_n) \one\{a_{i, j}(\PP_n) \ge 1 - b_{i,j}^{(M)}\}] + \E[a_{i,j}(\PP_n) \one\{a_{i, j}(\PP_n) < 1 - b_{i,j}^{(M)}\}]  \\
&\le \P(a_{i, j}(\PP_n) \ge 1 - b_{i,j}^{(M)}) + (1 - b_{i,j}^{(M)}) \P(a_{i, j}(\PP_n) < 1 - b_{i,j}^{(M)}) = 1 - b_{i,j}^{(M)} \P(a_{i, j}(\PP_n) < 1 - b_{i,j}^{(M)})
\end{align*}
and therefore,
\begin{align*}
\P(a_{i, j}(\PP_n) < 1 - b_{i,j}^{(M)}) &\le \P(Q_n^{(i)} \text{ not }(\bar{\mc S}_n^{(j)}(\PP_n),M)\text{-bounded}) / b_{i,j}^{(M)} \\
&\le \P(Q_n^{(i)} \text{ not }(\bar{\mc S}_n^{(j)}(\PP_n),M)\text{-bounded}) e^{M/2}.
\end{align*}
Whether $Q_n^{(i)}$ is $(\bar{\mc S}_n^{(j)}(\PP_n),M)$-bounded depends only on the configurations in boxes $Q_n^{(s)}$ for $s\in\NN_+(i)$. For each of them, $Q_n^{(s)} \cap \bar{\mc S}_n^{(j)}(\PP_n) \in\{Q_n^{(s)}\cap\PP_n, Q_n^{(s)}\cap\PP'_n\}$, i.e., there are less than $2^{\#\NN_+(i)}\le 2^{3^d}$ possibilities. With the union bound, this leads to
\begin{equation}\label{inequality_poisson_union}
\P(Q_n^{(i)}\text{ is not }(\bar{\mc S}_n^{(j)}(\PP_n),M)\text{-bounded}) \le 2^{3^d} \P(Q_n^{(i)}\text{ is not }(\PP_n,M)\text{-bounded}).
\end{equation}
From here, we can continue by using Markov's inequality and Mecke's formula. To simplify the notation we set $m_n := M+a_n + s_0$ and get
\begin{align}\label{inequality_probability_bad_box}
\begin{split}
	&\P(Q_n^{(i)}\text{ is not }(\PP_n,M)\text{-bounded}) = \P\Big(\min_{X\in Q_n^{(i)}\cap\PP_n} \PP_n\big(B_{(\frac{m_n}{n\kappa_d})^{1/d}}(X)\big) \le k \Big) \\
	&\le \E\bigg[\sum_{X\in Q_n^{(i)}\cap\PP_n} \one\big\{\PP_n\big(B_{(\frac{m_n}{n\kappa_d})^{1/d}}(X)\big) \le k\big\}\bigg] = n \int_{Q_n^{(i)}} \E\big[\one\big\{\PP_n\big(B_{(\frac{m_n}{n\kappa_d})^{1/d}}(x)\big) \le k-1\big\}\big] \d x \\
	&= n |Q_n^{(i)}| \sum_{i=0}^{k-1} \frac{m_n^i}{i!} e^{-m_n} \le \rho_n^\de |Q_n^{(i)}| k(1+M/a_n+s_0/a_n)^{k-1} e^{-M-s_0} \le 2k e^{-M-s_0}
\end{split}
\end{align}
for large enough $n$. With this, for the $j$th box of the arbitrary ordering, $Q_n^{(j)}$, we compute that
\begin{align*}
\P(Q^{(j)} \text{ is } (\PP_n,M)\text{-bad}) &= \P\Big(\bigcup_{i\in\NN_+(j)} \{a_{i,j}(\PP_n) < 1 - b_{i,j}^{(M)}\}\Big) \le \sum_{i\in\NN_+(j)} \P(a_{i,j}(\PP_n) < 1 - b_{i,j}^{(M)}) \\
&\le 3^d 2^{3^d+1} e^{M/2} ke^{-M-s_0} = 3^d 2^{3^d+1} k e^{-M/2-s_0},
\end{align*}
and thus, choosing $c:=3^d 2^{3^d+1} k$ suffices for the claim to hold.
\enp

\bep[Proof of Lemma \ref{lemma_prob_new_edges}]
First, we recall the events $E_j^\text{good}(\eta)$, $E_j^\text{good}(\eta)$ and $E_j^b$ from \eqref{equation_good_condition1} and \eqref{equation_good_condition2} for $\eta$ equal to either $\PP_n$ or $\PP'_n$.
Then, as a first step, we point out that by the tower property
\begin{align}\label{equality_recursive1}
\begin{split}
&\P\bigg(\bigcap_{j=1}^{\rho_n^\de} E_j^\text{good}(\PP_n)\cup \big(E_j^\text{bad}(\PP_n)\cap E_j^\text{good}(\PP'_n) \cap E_j^b\big) \biggm\vert \PP_n\bigg) \\
&= \E\bigg[\Big(\one_{E_{\rho_n^\de}^\text{good}(\PP_n)} + \one_{E_{\rho_n^\de}^\text{bad}(\PP_n)} \E\big[ \one_{E_{\rho_n^\de}^\text{good}(\PP'_n) \cap E_{\rho_n^\de}^b} \bigm\vert \PP_n, (\PP'_n)^{(\rho_n^\de-1)}\big]\Big) \\
&\qquad\ \prod_{j=1}^{\rho_n^\de-1} \big(\one_{E_j^\text{good}(\PP_n)} + \one_{E_j^\text{bad}(\PP_n)} \one_{E_j^\text{good}(\PP'_n) \cap E_j^b}\big) \biggm\vert \PP_n\bigg].
\end{split}
\end{align}
This gives an indication of the recursive approach to this proof. We start by working towards a bound of the inner conditional expectation after the equals sign of \eqref{equality_recursive1}. Fixing an arbitrary $j\in\{1,\dots,\rho_n^\de\}$, note that the $(\PP'_n,M)$-goodness of $Q_n^{(j)}$ does not depend on $\cup_{s\ge j} Q_n^{(s)}\cap\PP_n$ and therefore
\begin{equation}\label{equality_recursion2}
\P\big(E_j^\text{good}(\PP'_n), E_j^b \mid \PP_n, (\PP'_n)^{(j-1)}\big) = \P\big(E_j^\text{good}(\PP'_n), E_j^b \mid \PP_n^{(j-1)}, (\PP'_n)^{(j-1)}\big).
\end{equation}
Now, we can use the definition of goodness to arrive at
\begin{align}\label{inequality_recursion3}
\begin{split}
&\P(E_j^\text{good}(\PP'_n), E_j^b \mid \PP_n^{(j-1)}, (\PP'_n)^{(j-1)}) \\
&= \P\big(\cap_{i\in\NN_+(j)} \{a_{i,j}(\PP'_n) \ge 1-b_{i,j}^{(M)}\} \cap E_j^b \mid \PP_n^{(j-1)}, (\PP'_n)^{(j-1)}\big) \\
&\ge 1 - \sum_{i\in\NN_+(j)} \big(1-\P\big(a_{i,j}(\PP'_n) \ge 1-b_{i,j}^{(M)}, E_j^b \mid \PP_n^{(j-1)}, (\PP'_n)^{(j-1)}\big)\big).
\end{split}
\end{align}
Subsequently, the key step is to show that under $\cap_{s\le j-1} \big(E_s^\text{good}(\PP_n)\cup (E_s^\text{bad}(\PP_n)\cap E_s^\text{good}(\PP'_n))\big)$ for sufficiently large $n$
\begin{equation}\label{inequality_assertion_dense_proof}
\P\big(a_{i,j}(\PP'_n) \ge 1-b_{i,j}^{(M)}, E_j^b \mid \PP_n^{(j-1)}, (\PP'_n)^{(j-1)}\big) \ge 1 - 2 k e^{|s_0|} (e^{-M_0} - e^{-M/2^{4^d}}).
\end{equation}
Once \eqref{inequality_assertion_dense_proof} is established, we conclude the proof as follows. Continuing at \eqref{equality_recursion2} and \eqref{inequality_recursion3}, using that $\#\NN_+(j) = 3^d$, yields that
$$
\P\big(E_j^\text{good}(\PP'_n), E_j^b \mid \PP_n^{(j-1)}, (\PP'_n)^{(j-1)}\big) \ge 1 - 3^d 2 k e^{|s_0|} (e^{-M_0} - e^{-M/2^{4^d}}) = q_{M_0,M}.
$$
This lets us proceed at \eqref{equality_recursive1} to arrive at
\begin{align*}
&\P\bigg(\bigcap_{j=1}^{\rho_n^\de} \Big(E_j^\text{good}(\PP_n)\cup \big(E_j^\text{bad}(\PP_n)\cap E_j^\text{good}(\PP'_n) \cap E_j^b\big)\Big) \biggm\vert \PP_n\bigg) \\
&\ge \E\bigg[\Big(\one_{E_{\rho_n^\de}^\text{good}(\PP_n)} + \one_{E_{\rho_n^\de}^\text{bad}(\PP_n)} q_{M_0,M}\Big) \prod_{j=1}^{\rho_n^\de-1} \big(\one_{E_j^\text{good}(\PP_n)} + \one_{E_j^\text{bad}(\PP_n)} \one_{E_j^\text{good}(\PP'_n) \cap E_j^b}\big) \biggm\vert \PP_n\bigg] \\
&\ge q_{M_0,M} \E\bigg[ \prod_{j=1}^{\rho_n^\de-1} \big(\one_{E_j^\text{good}(\PP_n)} + \one_{E_j^\text{bad}(\PP_n)} \one_{E_j^\text{good}(\PP'_n) \cap E_j^b}\big) \biggm\vert \PP_n\bigg] \ge q_{M_0,M}^{\rho_n^\de},
\end{align*}
where the last inequality follows from repeating the previous steps $\rho_n^\de$ times.

It remains to prove the assertion stated in \eqref{inequality_assertion_dense_proof}. In order to do so, let $i\in\NN_+(j)$ be fixed. If $\{1,\dots,j-1\}\cap\NN_+(i)\neq\emptyset$, we can denote the largest index of an adjacent box of the box $Q_n^{(i)}$ that comes before $j$ in the ordering by $j_0 := \max(\{1,\dots,j-1\}\cap\NN_+(i))$. Note that $Q_n^{(j_0)}\cap(\bar{\mc S}_n^{(j)}(\PP'_n)$ can either be equal to $Q_n^{(j_0)}\cap\PP_n$ or $Q_n^{(j_0)}\cap\PP'_n$, resulting in two options that we can include in a similar way as was done in \eqref{inequality_poisson_union}. Then, we have that under $\cap_{s\le j-1} \big(E_s^\text{good}(\PP_n)\cup (E_s^\text{bad}(\PP_n)\cap E_s^\text{good}(\PP'_n))\big)$
\begin{align}\label{inequality_proof_with_recursion1}
\begin{split}
&\P(Q^{(i)}\text{ is $(\bar{\mc S}_n^{(j)}(\PP'_n),M)$-bounded}, E_j^b \mid \PP_n^{(j-1)}, (\PP'_n)^{(j-1)}) \\
&\ge \P(E_j^b \mid \PP_n^{(j-1)}, (\PP'_n)^{(j-1)}) - \P(Q^{(i)}\text{ is not }(\bar{\mc S}_n^{(j)}(\PP'_n),M)\text{-bounded} \mid \PP_n^{(j-1)}, (\PP'_n)^{(j-1)}) \\
&= \P(E_j^b) - \P(Q^{(i)}\text{ is not }(\bar{\mc S}_n^{(j)}(\PP'_n),M)\text{-bounded} \mid \PP_n^{(j_0)}, (\PP'_n)^{(j_0)}) \\
&\ge \P(E_j^b) - 2\P(Q^{(i)}\text{ is not }(\bar{\mc S}_n^{(j_0)}(\PP'_n),M)\text{-bounded} \mid \PP_n^{(j_0)}, (\PP'_n)^{(j_0)}) \\
&\ge \P(E_j^b) - 2 b_{i,j_0}^{(M)} = \P(E_j^b) - 2 b_{i,j-1}^{(M)}.
\end{split}
\end{align}
In the other case, i.e., if $\{1,\dots,j-1\}\cap\NN_+(i)=\emptyset$, we get
\begin{align}\label{inequality_proof_with_recursion2}
\begin{split}
&\P(Q^{(i)}\text{ is $(\bar{\mc S}_n^{(j)}(\PP'_n),M)$-bounded}, E_j^b \mid \PP_n^{(j-1)}, (\PP'_n)^{(j-1)}) \\
&= \P(Q^{(i)}\text{ is $(\PP'_n,M)$-bounded}, E_j^b) \ge \P(E_j^b) - 2ke^{-M-s_0}
\end{split}
\end{align}
for large $n$, where the last inequality follows from \eqref{inequality_probability_bad_box}. For completeness, note that we viewed $\PP_n^{(0)}$ and $(\PP'_n)^{(0)}$ as $\emptyset$. Additionally, with the tower property, it follows that
\begin{align*}
&\P(Q^{(i)}\text{ is $(\bar{\mc S}_n^{(j)}(\PP'_n),M)$-bounded}, E_j^b \mid \PP_n^{(j-1)}, (\PP'_n)^{(j-1)}) \\
&= \E\big[\P(Q^{(i)}\text{ is $(\bar{\mc S}_n^{(j)}(\PP'_n),M)$-bounded} \mid \PP_n^{(j)}, (\PP'_n)^{(j)}) \one\{E_j^b\} \bigm\vert \PP_n^{(j-1)}, (\PP'_n)^{(j-1)}\big] \\
&=\E\big[\P(Q^{(i)}\text{ is $(\bar{\mc S}_n^{(j)}(\PP'_n),M)$-bounded} \mid \PP_n^{(j)}, (\PP'_n)^{(j)}) \one\{E_j^b\} \\
&\qquad\ (\one\{a_{i,j}(\PP'_n) \ge 1 - b_{i,j}^{(M)}\} + \one\{a_{i,j}(\PP'_n) < 1 - b_{i,j}^{(M)}\}) \bigm\vert \PP_n^{(j-1)}, (\PP'_n)^{(j-1)}\big] \\
&\le\E\big[\one\{E_j^b\} \big(\one\{a_{i,j}(\PP'_n) \ge 1-b_{i,j}^{(M)}\} + (1 - b_{i,j}^{(M)}) \one\{a_{i,j}(\PP'_n) < 1 - b_{i,j}^{(M)}\}\big) \bigm\vert \PP_n^{(j-1)}, (\PP'_n)^{(j-1)}\big] \\
&\le (1 - b_{i,j}^{(M)}) \P(E_j^b) + b_{i,j}^{(M)} \P(a_{i,j}(\PP'_n) \ge 1 - b_{i,j}^{(M)}, E_j^b \mid \PP_n^{(j-1)}, (\PP'_n)^{(j-1)}).
\end{align*}
Note that similar to \eqref{inequality_probability_bad_box}, we can also show that $\P(E_j^b)\ge 1-2ke^{-M_0-s_0}$. Using this,  \eqref{inequality_proof_with_recursion1} and \eqref{inequality_proof_with_recursion2} as well as the definition of $b_{i,j}^{(M)}$ from \eqref{equation_probability_Mgood}, we arrive at
\begin{align*}
&\P(a_{i,j}(\PP'_n) \ge 1 - b_{i,j}^{(M)}, E_j^b \mid \PP_n^{(j-1)}, (\PP'_n)^{(j-1)}) \\
&\ge \frac{\P(E_j^b) - \max\{2 b_{i,j-1}^{(M)}, 2 k e^{-M-s_0}\} - (1 - b_{i,j}^{(M)}) \P(E_j^b)}{b_{i,j}^{(M)}} \ge \frac{b_{i,j}^{(M)} \P(E_j^b) - 2 k b_{i,j-1}^{(M)} e^{|s_0|}}{b_{i,j}^{(M)}} \\
&\ge \frac{b_{i,j}^{(M)}(1-2 k e^{-M_0-s_0}) - 2 k b_{i,j-1}^{(M)} e^{|s_0|}}{b_{i,j}^{(M)}} = 1-2 k e^{-M_0-s_0} - 2 k b_{i,j-1}^{(M)} e^{|s_0|} / b_{i,j}^{(M)} \\
&= 1 - 2 k e^{-M_0-s_0} - 2 k e^{-M (2^{-1-(\#\{s\in\NN_+(i) \colon s\le j\}-1)} - 2^{-1-\#\{s\in\NN_+(i) \colon s\le j\}})} e^{|s_0|} \\
&= 1 - 2 k e^{-M_0-s_0} - 2 k b_{i,j}^{(M)} e^{|s_0|} \ge 1 - 2 k e^{-M_0-s_0} - 2 k e^{-M 2^{-1-3^d}} e^{|s_0|} \\
&\ge 1 - 2 k e^{|s_0|} (e^{-M_0} - e^{-M/2^{4^d}}).
\end{align*}
\enp

\bep[Proof of Lemma \ref{lemma_lower_bound_good_event}]
For the first part, note that given $E_n^*$ the events $\{a_{i,\rho_n^\de}(\PP_n'') \ge 1-b_{i,\rho_n^\de}^{(M)}\}$ occur for all $i\in\mc N_+(\rho_n^\de)$. Thus, for sufficiently large $M$ and every $i\in\mc N_+(\rho_n^\de)$
\begin{align}\label{inequality_recursion_dense}
\begin{split}
0 &< 1 - b_{i,\rho_n^\de}^{(M)} \le a_{i, \rho_n^\de}(\PP''_n) = \P\big(Q_n^{(i)}\text{ is $(\bar{\mc S}_n^{(\rho_n^\de)}(\PP''_n),M)$-bounded}\bigm\vert \PP_n^{(\rho_n^\de)}, (\PP'_n)^{(\rho_n^\de)}\big) \\
&= \one\{Q^{(i)}\text{ is $(\bar{\mc S}_n^{(\rho_n^\de)}(\PP''_n),M)$-bounded}\} = \one\{Q_n^{(i)}\text{ is $(\PP''_n,M)$-bounded}\}
\end{split}
\end{align}
by measurability with respect to $\PP_n^{(\rho_n^\de)}, (\PP'_n)^{(\rho_n^\de)}$ and therefore, $Q_n^{(i)}$ is $(\PP''_n,M)$-bounded. Next, repeating this argument, it follows that $Q_n^{(i')}$ is $(\PP''_n,M)$-bounded for all $i'\in\mc N_+(\rho_n^\de-1) \setminus \mc N_+(\rho_n^\de)$. Note that $i'\not\in\mc N_+(\rho_n^\de)$ is an important requirement to be able to replicate the last two equalities in \eqref{inequality_recursion_dense} in this case. Afterwards, we consider $i''\in\mc N_+(\rho_n^\de-2) \setminus \big(\mc N_+(\rho_n^\de-1) \cup \mc N_+(\rho_n^\de-1)\big)$. We can repeat this until all boxes have been dealt with and we conclude the first part of the proof of Lemma \ref{lemma_lower_bound_good_event} by deducing from this that given $E_n^*$, the event that $Q_n^{(j)}$ is $(\PP''_n,M)$-bounded holds for all $j\in\{1,\dots,\rho_n^\de\}$.

For the second part, the tower property yields
\begin{align}\label{inequality_bernoulli_separation}
\begin{split}
\P(E_n^* \mid \PP_n) &= \E\bigg[\prod_{j=1}^{\rho_n^\de} \big(\one_{\{E_j^\text{good}(\PP_n)\}} \one_{\{X_{Q_n^{(j)},\e} = 0\}} + \one_{\{E_j^\text{bad}(\PP_n)\}} \one_{\{E_j^\text{good}(\PP'_n)\cap E_j^b\}} \one_{\{X_{Q_n^{(j)},\e} = 1\}}\big) \biggm\vert \PP_n\bigg] \\
&= \E\bigg[\prod_{j=1}^{\rho_n^\de} \E\Big[\one_{\{E_j^\text{good}(\PP_n)\}} \one_{\{X_{Q_n^{(j)},\e} = 0\}} \\
&\qquad\qquad\qquad + \one_{\{E_j^\text{bad}(\PP_n)\}} \one_{\{E_j^\text{good}(\PP'_n)\cap E_j^b\}} \one_{\{X_{Q_n^{(j)},\e} = 1\}} \Bigm\vert \PP_n, \PP'_n\Big] \biggm\vert \PP_n\bigg].
\end{split}
\end{align}
From this point, using the independence of $(X_{Q_n^{(j)},\e})_j$ of all Poisson point processes and the measurability of $E_j^\text{good}(\PP_n)$, $E_j^\text{good}(\PP'_n)$ and $E_j^b$ with respect to $\sigma(\PP_n,\PP'_n)$, we can compute that it almost surely holds that
\begin{align*}
&\E\Big[\one_{\{E_j^\text{good}(\PP_n)\}} \one_{\{X_{Q_n^{(j)},\e} = 0\}} + \one_{\{E_j^\text{bad}(\PP_n)\}} \one_{\{E_j^\text{good}(\PP'_n)\cap E_j^b\}} \one_{\{X_{Q_n^{(j)},\e} = 1\}} \Bigm\vert \PP_n, \PP'_n\Big] \\
&= \e^{\one_{\{E_j^\text{good}(\PP_n)\}}}  (1-\e)^{\one_{\{E_j^\text{bad}(\PP_n)\}} \one_{\{E_j^\text{good}(\PP'_n)\cap E_j^b\}}} (\one_{\{E_j^\text{good}(\PP_n)\}} + \one_{\{E_j^\text{bad}(\PP_n)\}} \one_{\{E_j^\text{good}(\PP'_n)\cap E_j^b\}}) \\
&\ge \e^{\one_{\{E_j^\text{good}(\PP_n)\}}} (1-\e) (\one_{\{E_j^\text{good}(\PP_n)\}} + \one_{\{E_j^\text{bad}(\PP_n)\}} \one_{\{E_j^\text{good}(\PP'_n)\cap E_j^b\}}).
\end{align*}
Next, we revisit \eqref{inequality_bernoulli_separation} and continue with
\begin{align*}
&\P(E_n^* \mid \PP_n)\\
 &\ge \E\bigg[\e^{J_n^M}  (1-\e)^{\rho_n^\de} \prod_{j=1}^{\rho_n^\de} \big(\one_{\{E_j^\text{good}(\PP_n)\}} + \one_{\{E_j^\text{bad}(\PP_n)\}} \one_{\{E_j^\text{good}(\PP'_n)\cap E_j^b\}} \big) \biggm\vert \PP_n\bigg] \\
&\ge \E\bigg[\e^{\delta\rho_n^\de} (1-\e)^{\rho_n^\de} \one_{\{J_n^M < \delta\rho_n^\de\}} \prod_{j=1}^{\rho_n^\de} \big(\one_{\{E_j^\text{good}(\PP_n)\}} + \one_{\{E_j^\text{bad}(\PP_n)\}} \one_{\{E_j^\text{good}(\PP'_n)\cap E_j^b\}} \big) \biggm\vert \PP_n\bigg] \\
&\ge \e^{\delta\rho_n^\de} (1-\e)^{\rho_n^\de}  \bigg(\E\bigg[\prod_{j=1}^{\rho_n^\de} \big(\one_{\{E_j^\text{good}(\PP_n)\}} + \one_{\{E_j^\text{bad}(\PP_n)\}} \one_{\{E_j^\text{good}(\PP'_n)\cap E_j^b\}} \big) \biggm\vert \PP_n\bigg] \\
&\qquad\qquad\qquad\qquad\quad - \P(J_n^M \ge \delta\rho_n^\de \mid \PP_n)\bigg).
\end{align*}
Here, Lemma \ref{lemma_prob_new_edges} yields
\begin{equation}\label{inequality_good_event2}
\P(E_n^* \mid \PP_n) \ge \e^{\delta\rho_n^\de} (1-\e)^{\rho_n^\de}  \big(q_{M_0,M}^{\rho_n^\de} - \P(J_n^M \ge \delta\rho_n^\de \mid \PP_n)\big).
\end{equation}
\enp

\bep[Proof of Lemma \ref{lemma_prob_resampling_between_boxes}]
For a box $Q\in\mathcal{Q}_n$, we divide $\partial_n Q$ into a grid consisting of boxes with side length $u_n := (\frac{a_n+s_0}{n\kappa_d})^{1/d}$ and call this collection of boxes $\mathcal{W}_n(\partial_n Q)$. We denote the total collection of these boxes by $\overline{\mathcal{W}}_n := \cup_{Q\in\mathcal{Q}_n} \mathcal{W}_n(\partial_n Q)$. The volume of $\partial_n Q$ for $Q\in\mathcal{Q}_n$ can be bounded by $5 d (\rho_n^\de)^{-(d-1)/d} t_n$
for large $n$ and therefore,
\begin{equation}\label{limit_zero_boxes_partial}
|\partial_n Q|/ |Q| \le 5 d (\rho_n^\de)^{-(d-1)/d} t_n \rho_n^\de = 5 d (\rho_n^\de)^{1/d} t_n = e^{-a_n/d} a_n^{(k-1)/d} (a_n+w_n)^{1/d} / \kappa_d^{1/d} \overset{n\to\infty}{\longrightarrow} 0.
\end{equation}

Next, we can proceed with the same strategy that was previously employed to prove Lemma \ref{lemma_probability_bad_boxes}. We use the label set $\mathcal{L} := \{1,2,3\}^d$ to achieve that between two boxes of the same label, there are always two boxes labeled differently, where for simplicity, we assume that the number of boxes along each axis is divisible by $3$. We reuse the notation $\overline{\mathcal{W}}_n^{(l)}$ for the boxes of label $l\in\mathcal{L}$. Let $M,\delta>0$. Now, we can union over all labels and combine this with the union bound to arrive at
\begin{align*}
	&\P\bigg(\sum_{ X \in \PP_n  \cap (\cup_{Q\in\mathcal{Q}_n} \partial_n Q) } M \wedge \xi_n(X, \PP_n) \ge \delta\rho_n^\de\bigg) \\
	&\le \sum_{l\in\mathcal{L}} \P\bigg(\sum_{ X \in \PP_n  \cap (\cup_{W\in\overline{\mathcal{W}}_n^{(l)}} W) } M \wedge \xi_n(X, \PP_n) \ge \delta\rho_n^\de/3^d\bigg).
\end{align*}
For each $W\in\overline{\mathcal{W}}_n^{(l)}$, we assert that the maximal number of Poisson points in $X\in W\cap\PP_n$ with $R_k(X,\PP_n) \ge u_n$ is bounded by some $c := c(d,k)>0$. This follows similarly as in \eqref{inequality_nodes_nearest_neighbor_large_stabilization}. We go through nodes $W\cap\PP_n$ one by one and label some of them in the same manner as in Section \ref{section_applications_thermodynamic_kNN}. The only difference is that we can only argue that a fraction of $1/2^d$ of the volume of each constructed disjoint ball is in $W$, to account for vertices close to the boundary of $W$. This means the bound is computed by
$$\frac{k|W|}{\kappa_d (u_n/2)^d / 2^d} = k4^d/\kappa_d =: c.$$
Using this, for each $l\in\mathcal{L}$, we compute
\begin{align*}
&\P\bigg(\sum_{ X \in \PP_n  \cap (\cup_{W\in\overline{\mathcal{W}}_n^{(l)}} W) } M \wedge \xi_n(X, \PP_n) \ge \delta\rho_n^\de/3^d\bigg) \\
&\le \P\bigg(M \sum_{ X \in \PP_n  \cap (\cup_{W\in\overline{\mathcal{W}}_n^{(l)}} W) } \one\{R_k(X, \PP_n) \ge u_n\} \ge \delta\rho_n^\de/3^d\bigg) \\
&\le \P\bigg( \sum_{W\in\overline{\mathcal{W}}_n^{(l)}} \one\{\max_{X\in W\cap\PP_n} R_k(X, \PP_n) \ge u_n\} \ge \delta\rho_n^\de/(3^d c M)\bigg).
\end{align*}
With the goal of using the spatial independence to invoke a binomial concentration bound, we combine Markov's inequality and Mecke's equation, which yields for each $W\in\overline{\mathcal{W}}_n^{(l)}$ and $n$ large
\begin{align*}
\P\big(\max_{X\in W\cap\PP_n} R_k(X, \PP_n) \ge u_n\big) &\le n \int_W \P(R_k(x, \PP_n\cup\{x\}) \ge u_n) \d x = n \int_W \P(\PP_n(B_{u_n}(x)) < k) \d x \\
&= n\int_W \sum_{i=0}^{k-1} e^{- n u_n^d \kappa_d} \frac{(n u_n^d \kappa_d)^i}{i!} \d x \le n\int_W k e^{- (a_n + s_0)} (a_n + s_0)^{k-1} \d x \\
&= n|W| k e^{- (a_n + s_0)} (a_n + s_0)^{k-1} = |W| \rho_n^\de k e^{-s_0} (1 + s_0/a_n)^{k-1}.
\end{align*}
With the binomial bound from \cite[Lemma 1.2]{poisson_conc} and the computations from \eqref{limit_zero_boxes_partial} we arrive at
\begin{align*}
&\P\Big( \sum_{W\in\overline{\mathcal{W}}_n^{(l)}} \one\{\max_{X\in W\cap\PP_n} R_k(X, \PP_n) \ge u_n \} \ge \delta\rho_n^\de/(3^d c M)\Big) \\
&\le \exp\bigg(- \frac{\delta\rho_n^\de}{3^d 2 c M} \log\Big(\frac{\delta\rho_n^\de/(3^d c M)}{|W| \rho_n^\de k e^{-s_0} (1 + s_0/a_n)^{k-1} \rho_n^\de |\partial_n Q_n^{(1)}|/|W| }\Big)\bigg) \\
&= \exp\bigg(- \frac{\delta\rho_n^\de}{3^d 2 c M} \log\Big(\frac{\delta / (3^d c M)}{k e^{-s_0} (1 + s_0/a_n)^{k-1} \underbrace{|\partial_n Q_n^{(1)}| / |Q_n^{(1)}|}_{\overset{n\to\infty}{\longrightarrow}0}}\Big)\bigg).
\end{align*}
Thus, $\limsup_{n\to\infty} \frac1{\rho_n^\de} \log\P(H_{n,M}^{\text{err},\partial}(\PP''_n) \ge \delta) = -\infty$.

For the second part, we proceed roughly in the same fashion. But first, we note that for additionally $M_0,\tilde\delta>0$,
\begin{align*}
&\P(H_{n,M,M_0}^{\text{err}, \mc J}(\PP_n, \PP'_n, \PP''_n) \ge \delta) \\
&\le \P(H_{n,M,M_0}^{\text{err}, \mc J}(\PP_n, \PP'_n, \PP''_n) \ge \delta, J_n^M < \tilde\delta \rho_n^\de) + \P(J_n^M \ge \tilde\delta \rho_n^\de).
\end{align*}
In the following computations, we will use the upper bound for the binomial coefficient $\binom{a}{b} \le (ea/b)^b$, see \cite[Section 1.2.6 Exercise 67]{knuth1997}, Applied here, it yields
$$\binom{\rho_n^\de}{\tilde\delta\rho_n^\de} \le (e\rho_n^\de / (\tilde\delta \rho_n^\de))^{\tilde\delta \rho_n^\de} = (e / \tilde\delta)^{\tilde\delta \rho_n^\de},$$
where we assume that $\tilde\delta\le 1/2$ and that the pair of numbers occurring in the binomial coefficient are both positive integers. Now, we continue with
\begin{align*}
&\P(H_{n,M,M_0}^{\text{err}, \mc J}(\PP_n, \PP'_n, \PP''_n) \ge \delta, J_n^M < \tilde\delta \rho_n^\de) \\
&\le \P\bigg(\bigcup_{\mc A\subseteq \mc Q_n, \#\mc A < \tilde\delta \rho_n^\de} \bigg\{\sum_{ X \in \PP_n''  \cap (\cup_{Q\in\mc A} Q ) } M_0 \wedge \xi_n(X, \PP_n'') \ge \delta\rho_{n, k}^\de\bigg\}\bigg) \\
&\le \sum_{i=1}^{\tilde\delta \rho_n^\de} \sum_{\mc A\subseteq \mc Q_n, \#\mc A = i} \P\bigg(\sum_{ X \in \PP_n''  \cap (\cup_{Q\in\mc A} Q ) } M_0 \wedge \xi_n(X, \PP_n'') \ge \delta\rho_{n, k}^\de\bigg) \\
&\le \tilde\delta \rho_n^\de (e / \tilde\delta )^{\tilde\delta \rho_n^\de} \P\bigg(\sum_{ X \in \PP_n''  \cap (\cup_{i=1}^{\tilde\delta \rho_n^\de} Q_n^{(i)} ) } M_0 \wedge \xi_n(X, \PP_n'') \ge \delta\rho_{n, k}^\de\bigg).
\end{align*}
Next, we cover $\cup_{i=1}^{\tilde\delta \rho_n^\de} Q_n^{(i)}$ with cubes of side length $u_n$ and consistently with prior convention, denote this collection by $\overline{W}_n$. Then, $\#\overline{W}_n = \tilde\delta \rho_n^\de |Q_n^{(1)}| u_n^{-d} = \tilde\delta u_n^{-d}$, which we assume to be an integer. Next, we can simply introduce the same labeling as for the first part of this proof and by the same calculations as in the first part, we get
\begin{align*}
	&\P\bigg(\sum_{ X \in \PP_n''  \cap (\cup_{i=1}^{\tilde\delta \rho_n^\de} Q_n^{(i)} ) } M_0 \wedge \xi_n(X,\PP''_n) \ge \delta\rho_{n, k}^\de\bigg) \\
	&\le \sum_{l\in\mathcal{L}} \exp\Big(- \frac{\delta\rho_n^\de}{3^d 2 c M_0} \log\Big(\frac{\delta\rho_n^\de/(3^d c M_0)}{|W| \rho_n^\de k e^{-s_0} (1 + s_0/a_n)^{k-1} \tilde\delta u_n^{-d}}\Big)\Big) \\
	&\le 3^d \exp\Big(- \frac{\delta\rho_n^\de}{3^d 2 c M_0} \log\Big(\frac{\delta}{\tilde\delta 3^d c M_0 k e^{- s_0} (1 + s_0/a_n)^{k-1}}\Big)\Big)
\end{align*}
for large $n$. When choosing $\tilde\delta=\delta/\log M$, this yields
\begin{align*}
&\frac1{\rho_n^\de} \log \P(H_{n,M,M_0}^{\text{err}, \mc J}(\PP_n, \PP'_n, \PP''_n) \ge \delta, J_n^M < \tilde\delta \rho_n^\de) \\
&\le \frac1{\rho_n^\de}\log\Big(\frac{\delta \rho_n^\de 3^d}{ \log M}\Big) + \frac{\delta\log(e (\log M) / \delta)}{\log M} - \frac{\delta}{3^d 2 c M_0} \log\Big(\frac{\log M}{ 3^d c M_0 k e^{- s_0} (1 + s_0/a_n)^{k-1}}\Big) \\
&\overset{n\to\infty}{\longrightarrow} \frac{\delta\log(e (\log M) / \delta)}{\log M} - \frac{\delta}{3^d 2 c M_0} \log\Big(\frac{\log M}{ 3^d c M_0 k e^{- s_0} }\Big) \overset{M\to\infty}{\longrightarrow} -\infty.
\end{align*}
On the other hand, from Lemma \ref{lemma_probability_bad_boxes}, for $\tilde\delta=\delta/\log M$, we can deduce that
\begin{align*}
&\frac1{\rho_n^\de} \log\P(J_n^M \ge \tilde\delta\rho_n^\de) \le -\frac{\delta/\log M}{5^d 2} \log\Big(\frac{\delta e^{M/2} / \log M}{15^d 2^{3^d}}\Big) \\
&= \frac{-M/2 + \log\log M}{\log M} \frac{\delta}{5^d 2} - \frac{\delta/\log M}{5^d 2} \log\Big(\frac{\delta}{15^d 2^{3^d}}\Big)  \overset{M\to\infty}{\longrightarrow} -\infty
\end{align*}
and conclude the assertion.
\enp

\subsection*{Acknowledgment.}
The authors thank T.\ Owada for very fruitful discussions about the lower large deviations in the sparse regime. Further, DW would like to acknowledge the financial support of the CogniGron research center and the Ubbo Emmius Funds (Univ.\ of Groningen). 

\bibliography{refs}
\bibliographystyle{abbrv}

\end{document}